\documentstyle[12pt]{article}
\def\Bbb{\bf}
\def\be{\begin{equation}}
\def\ee{\end{equation}}
\def\bea{\begin{eqnarray*}}
\def\eea{\end{eqnarray*}}

\newcommand{\mapsonto}{\raisebox{3pt}
                            {\begin{picture}(14,2)\put(0,1){\vector(1,0){10}}
                                              \put(0,1){\vector(1,0){13}}
                         \end{picture}}}
\newtheorem{thm}{Theorem}[section]
\newtheorem{propn}[thm]{Proposition}
\newtheorem{cor}[thm]{Corollary}
\newtheorem{conjecture}{Conjecture}
\newcommand{\tr}{\Lambda}
\newtheorem{lemma}[thm]{Lemma}
\newtheorem{defn}{Definition}
\newenvironment{example}{ \medskip \noindent {\bf Example.}}{ \bigskip }
\newenvironment{remark}{\medskip
\noindent {\bf Remark.}}{\hfill \rule{.5em}{1em} \\}
\newenvironment{proof}{\medskip \noindent {\bf Proof.}}{\hfill \rule{.5em}{1em}
\\}
\def\b{\overline}

\def\+{\oplus}

\def\*{^{\ast}}
\def\dt{\frac{d}{dt}}

\def\vol{~d\mu} 
\def\volm{d\mu} 

\def\s{\sigma}

\def\om{\omega}

\def\bp{{\Bbb P}}
\def\BP{{\Bbb CP}}

\def\br{{\Bbb R}}
\def\bc{{\Bbb C}}
\def\bcp{{\Bbb CP}}
\def\B{{\cal B}}

\def\O{{\cal O}}
\begin{document}
\sloppy
\title{Existence and Deformation Theory for  Scalar-Flat K\"ahler Metrics
on Compact Complex Surfaces}

\author{\parbox{2in}{\center Claude LeBrun\thanks{Supported
in part by  NSF grant DMS-9003263.
{ Running title:}
{\sc Scalar-Flat K\"ahler Surfaces}
}\\SUNY Stony
 Brook} ~~~
  {\em and}~~~
 \parbox{2in}{\center Michael Singer\\Lincoln College, Oxford} }

\date{}
\maketitle





\section{Introduction}
\subsection{Motivation}

The classical uniformization theorem  provides a complete
translation dictionary
for the etymologically unrelated  languages of
complex 1-manifolds and constant curvature Riemannian 2-manifolds.
In higher  dimensions, there are a number of natural
ways in which one might try to generalize this remarkable theorem;
unfortunately,  these various potential generalizations remain,
for the most part,
programs  rather than established bodies of fact.
However, the subject  of the present article, namely
{\sl the existence problem for  zero-scalar-curvature
K\"ahler metrics on compact complex 2-manifolds},
occupies the  cross-roads of several such avenues of  research; and by
  clearing up a substantial piece
 of this problem, we thereby hope to  facilitate the flow
of traffic heading on to
a number interesting destinations.

 Purely in the context of Riemannian geometry, the most
optimistic programs to
 generalize the classical uniformization
theorem would try to equip every compact smooth manifold of
a given  dimension
with a (small!) class of ``optimal'' or ``canonical'' metrics.
In dimension four, one of the most natural versions  would have us seek
 extrema (or perhaps just critical points) of the   squared $L^2$-norm
$${\cal R}(g)= \int_M\|R\|^2 \vol$$
of the Riemann curvature
tensor $R$
over the space of  smooth Riemannian metrics $g$ on a given
smooth, compact, oriented
4-manifold $M$. Using the Chern-Gauss-Bonnet
formulas for the Euler characteristic $\chi$ and signature $\tau$
of our manifold $M$, one may easily  show \cite{L5} that
\bea {\cal R}(g)&=& -8\pi^2(3\tau +\chi )+\int_M
(4\|W_{+}\|^2+\frac{s^2}{12})\vol
\\ &\geq& -8\pi^2(3\tau +\chi )  ~,\eea
with equality iff $W_+=s=0$; here $s$ denotes the scalar curvature
 and
$W_+$ the self-dual Weyl curvature (cf. \S \ref{asd}) of $g$.
Metrics with $W_+=s=0$,  when they exist, are thus  absolute minima of
${\cal R}$, and it is therefore natural to try
to determine which
manifolds $M$ can admit such metrics. However, if the intersection form
$$\cup : H^2(M, {\Bbb R})\times H^2(M, {\Bbb R})\to {\Bbb R}$$
is indefinite, an elementary  Weitzenb\"ock argument \cite{L0}
shows that such a manifold admits an integrable complex structure
with respect to which the metric is K\"ahler; conversely \cite{G}, any
K\"ahler manifold of complex dimension 2 with $s\equiv 0$ automatically
satisfies $W_+=0$ and has indefinite intersection form.
 Thus the problem of minimizing $\cal R$ on a smooth
manifold leads us quite naturally\footnote{It should be pointed
out that many manifolds which do not admit scalar-flat K\"ahler
metrics nonetheless admit metrics which are absolute minima
of $\cal R$. In particular \cite{besse},
any {\em Einstein metric} on a compact 4-manifold
provides an absolute
minimum of $\cal R$. It is this fact which explains much of the
current interest in this Riemannian functional.}
to the problem of classifying
compact  K\"ahler   manifolds of complex dimension 2 and
scalar curvature zero--- henceforth referred to as scalar-flat K\"ahler
surfaces.

A related  four-dimensional program would instead seek to optimize the
{\em conformal} geometry of Riemannian metrics by  seeking
to minimize the conformally-invariant squared $L^2$-norm
$${\cal W}(g)=\int_M\| W\|^2\vol $$
of the  conformal curvature over the space of conformal classes of
Riemannian
metrics on $M$. Since
$${\cal W}(g)= - 12\pi^2 \tau + 2\int_M\| W_+\|^2\vol  ~,$$
anti-self-dual metrics (i.e. metrics satisfying $W_+=0$)
 are obviously absolute minima of ${\cal W}$,
 scalar-flat  K\"ahler surfaces again provide examples of
absolute minima. Note that
 while there are strong topological  constraints on
anti-self-dual metrics with  non-negative scalar curvature
\cite{DF}  \cite{flo} \cite{L0} \cite{L} \cite{poon}, the situation
is radically different once the scalar curvature condition is
dropped; in fact,
it has recently been shown \cite{taubes} that
the  obstructions to the
existence of anti-self-dual metrics on any  oriented smooth 4-manifold
are so weak that they can always be killed off by ``blowing up points,''
in the differentiable
sense of taking connected sums with enough $\overline{\Bbb CP}_2$'s. Our own
results in this article will
have  something of a  similar ring to them--- while
there are a number of obstructions  to the existence of
scalar-flat K\"ahler metrics on a compact complex surface, we will
see that all but the crudest can be killed off
by blowing up points.

In search of a natural bridge between  complex and differential geometry,
Calabi has proposed the problem
 of representing K\"ahler classes
on  compact complex manifolds
 by  K\"ahler metrics of constant scalar
curvature. Here again there is a natural variational approach
to the problem,  since such metrics
are  absolute minima of the functional
$${\cal C}(g)=\int_M s^2 \vol$$
among metrics in a fixed K\"ahler class; more generally,
 critical points of this functional have come to be
known \cite{besse} \cite{cal}
as {\em extremal K\"ahler metrics}.   However,
the existence of constant
scalar curvature K\"ahler metrics is,
in general,
obstructed \cite{burnsbart}\cite{besse}\cite{fut}\cite{cal2},
and the known obstructions will necessarily play a central r\^ole
in the present article--- although  perhaps not quite in
 the way the reader might  expect.
It is  hoped that our
 present existence results
will  provide  a useful way station, {\sl en route} to a more general
understanding  of    Calabi's problem.

   From a quite  different perspective, namely
that  of Hawking's Euclideanization program in gravitational physics, a
fundamental problem is that of classifying compact Riemannian solutions of
the {\em Einstein-Maxwell equations}
$$-{\textstyle \frac{1}{2}}
 r^{\sharp}=\mbox{{\em trace-free part}~}(F^{\sharp}\circ F^{\sharp}) $$
       $$     dF = d\star F=0 $$
governing the interaction of the gravitational field, represented by
a Riemannian metric $g$, with
the electromagnetic field, represented by a harmonic 2-form  $F$;
here the metric has been used to identify the Ricci curvature and
electromagnetic field with endomorphisms $r^{\sharp}$  and  $F^{\sharp}$
of the tangent bundle.
Any scalar-flat  K\"ahler surface provides a solution of these
equations once one sets
$$F=\rho +{\textstyle \frac{1}{4}}\omega ~,$$
where $\omega$ and $\rho$ are respectively the K\"ahler and Ricci
forms\footnote{If $F$ is to be viewed as a ${\bf U}(1)$-gauge field,
as it must be in realistic physical theories,   our K\"ahler metric
must  also be of  Hodge type--- that is, cohomologous to
the  metric induced by some projective embedding.};
conversely \cite{fla}, these are essentially the only solutions with
$W_+=0$. Thus the classification problem for scalar-flat K\"ahler
surfaces may be seen as part of a quest to classify
electro-gravitational
instantons.

Finally, the Penrose twistor correspondence \cite{P} gives  quite a different
way of generalizing the conformal surface/complex curve dictionary
to  dimension four. If $M$ is a smooth oriented  4-manifold
and $$[g]=\{ e^u g\}$$
is a conformal class of Riemannian metrics on $M$, the space $Z$
of orthogonal complex structures on $TM$ compatible with the orientation
is  an almost-complex manifold  of real dimension 6; it then turns out
\cite{AHS} that
$Z$ is a complex manifold iff $(M, [g])$ satisfies $W_+=0$. In this case,
$Z$ is called the {\em twistor space} of $(M, [g])$, and it turns out that
both $M$ and its anti-self-dual conformal structure can be reconstructed from
this complex manifold.
In particular, every scalar-flat K\"ahler surface $(M,g, J)$ has associated to
it a compact complex 3-fold $Z$; moreover, $(M,J)$ is
naturally  a complex submanifold of  $Z$.
While this construction has elsewhere served \cite{LP}\cite{taubes}
 primarily as an excellent source of  complex
3-folds with various ``pathologically'' non-K\"ahlerian
properties,  the deformation theory of $Z$
 will here serve as our
guiding light
as we trek through the  realm  of scalar-flat K\"ahler geometry.

\subsection{Outline}
{\sl We now provide the reader with a statement of the central result
of the paper, followed by an indication of the structure of the argument. }
\bigskip

\noindent {\bf Main Theorem}~~
{\em Let $M$ be a compact complex surface which admits
a K\"ahler metric whose scalar curvature has integral zero. Suppose
$\pi_1 (M)$ does not contain an Abelian subgroup of finite index.
 Then if $M$ is blown up at sufficiently many points, the
resulting surface $\tilde{M}$ admits scalar-flat K\"ahler metrics.}

\begin{description}
\item{\S \ref{flag} } We study the
behavior of the solution space under small deformations
of complex structure of the complex surface in question.
Our approach uses the twistor correspondence and
a modified version of Kodaira-Spencer theory.
This deformation theory is generally obstructed, but
we are able to  describe the obstructions completely in terms of
the Futaki character of the algebra of holomorphic vector fields.
\item{\S\ref{foo} } We compute the Futaki character for all relevant
complex surfaces, and use this to show that the deformation theory is
unobstructed for all non-minimal surfaces.
 \item{\S\ref{next} } We describe a
 large class of exact solutions  previously found in \cite{L2}, and
an improvement on those results  which
gives a classification of scalar-flat  K\"ahler surfaces
with  non-trivial automorphism algebra.
 \item{\S\ref{tag} } Using the bimeromorphic classification theory
of surfaces, we show that every surface satisfying the
hypotheses of the Main Theorem has blow-ups which are
  arbitrarily small
deformations of surfaces on which we have exact solutions. Applying our
deformation theory then proves the Main Theorem.
\end{description}

\noindent  It should be emphasized  that
the  Main Theorem's fundamental-group  hypothesis
reflects
the limitations of our
techniques rather than a
known  obstruction
to the existence of scalar-flat K\"ahler metrics.
In fact, the article concludes with some
 speculations  (Conjecture \ref{better}, \S\ref{tag})
to the effect that this restriction is essentially
superfluous.

\medskip

\noindent Incidental to the main course of the argument,
we will also
encounter  a remarkable empirical relationship,
observed (\S\ref{next}) in two quite different classes of explicit examples,
which seems to link the existence
problem for scalar-flat K\"ahler metrics
to the  stability of vector bundles with parabolic structure
in the sense of Seshadri \cite{sesh}; it is our belief, as expressed in
Conjecture 2,   \S\ref{tag}, that this relationship will actually turn out to
hold in complete
generality.

\subsection{Notation and Conventions}
For the purposes of this
section, $(M,g)$ will denote
an oriented Riemannian $2m$-manifold, although in the sequel we will
specialize to the case $m=2$.
We  use
$\volm$ to denote the  volume form of $g$, and
$\nabla$ to denote its Levi-Civit\`{a} connection.
The  $C^{\infty}$ sections
of any smooth vector bundle ${\cal V}\to M$
will be denoted by ${\cal  E}({\cal V})$.
The operation of raising (lowering) an index will  be indicated by
$\sharp$ ($\flat$).
We give
 the  curvature tensor $R$ the usual Riemannian sign:
\begin{equation}
(\nabla_{c}\nabla_{d}-\nabla_{d}\nabla_{c})\xi^a={R^a}_{bcd}\xi^b
{}~~~\forall~\xi\in {\cal E}(TM).
\label{n1}
\end{equation}
The Ricci tensor ${R^c}_{acb}$
is denoted by $r_{ab}$ and the scalar curvature ${R^{ab}}_{ab}$ by $s$.

The pointwise inner product induced by $g$ on the tensor bundles will be
denoted
$(~,~)$, whereas the
 the global $L^2$ inner product will be denoted
by $\langle~,~\rangle$.  The formal
adjoint of $d$ with respect to this inner product will be denoted
$\delta$ and
is given by the usual formula
\begin{equation}
\delta=-\star d\star.\label{n2}
\end{equation}

The metric $g$ is said to be {\em K\"{a}hler} if its holonomy
group is conjugate
to a subgroup of ${\bf U}(m)\subset {\bf O }(2m)$.
More concretely, this means that there is a compatible
parallel almost-complex
structure $J$:
\begin{equation}
J^2=-{\bf 1}, ~~g(J\xi , J\eta )=g(\xi , \eta)~\forall \xi , \eta\in TM,
{}~~~\nabla J=0.\label{n3}
\end{equation}
If the holonomy of $g$ happens to be
 a proper subgroup of ${\bf U}(m)$,
there may be more than one $J$ which satisfies (\ref{n3}); nonetheless,
when we speak
of a K\"ahler metric we will henceforth
always assume that a {\em particular choice} of $J$ has been made.
We therefore have a decomposition
$ {\Bbb C}\otimes TM=T^{1,0}\oplus T^{0,1} $
into the $\pm i$ eigenspaces of $J$,
thereby inducing a decomposition
\be\wedge^r_{\Bbb C}=\bigoplus_{r=p+q}\wedge^{p,q}\label{decomp} \ee
of the bundle of $r$-forms into forms of type
$(p,q)$, as defined by $\wedge^{p,q}:= (\wedge^pT^{1,0})^{\ast}\otimes
(\wedge^qT^{0,1})^{\ast}$;
in particular, $J$ induces
a ``standard'' orientation of $M$ by requiring that the
$2m$-form $i^m\varphi \wedge\bar{\varphi}$
be positive for any non-zero element $\varphi$ of the {\em canonical
line bundle} $\kappa = \wedge^{m,0}$. For brevity, we will use
${\cal E}^r$ and ${\cal E}^{p,q}$ to respectively denote ${\cal E}(\wedge^r)$
and ${\cal E}(\wedge^{p,q})$.
Because $\nabla$ is torsion-free, $[ {\cal E} (T^{1,0}),{\cal E} (T^{1,0})]
\subset {\cal E} (T^{1,0})$, and the
 the Newlander-Nirenberg \cite{NN} theorem therefore
asserts that $J$ is integrable---
that is, there exists a system of
 local coordinates for which $J$ becomes the standard
almost-complex structure
on ${\Bbb C}^m$, making $M$ a complex $m$-manifold.
The K\"{a}hler form $\omega$ and Ricci form $\rho$ are then
defined by the formulae
\begin{equation}
\omega(X,Y)=g(J\xi, Y),~~~\rho(\xi,\eta)=r(J\xi, \eta)
{}~~~\forall \xi,\eta\in TM~;\label{n4}
\end{equation}
both are real closed forms of type (1,1).
Conversely, a closed real (1,1)-form on complex manifold is a K\"ahler
form iff the symmetric form $g$ it
defines implicitly via (\ref{n4}) is positive definite.
The  deRham class
$[\omega ]\in H^{2}(M)$ of the K\"ahler form
is called the K\"{a}hler class. It is a central fact
of K\"ahler geometry
that $\rho$ is exactly the curvature   of the  Chern connection
on $\kappa^{-1}=\wedge^mT^{1,0}$; in particular,
$\rho$ is completely determined by $J$ and $\volm$ alone,
and  the deRham class of $\rho /2\pi$
is just the first Chern class $c_{1}(M):=c_1(T^{1,0}M)=c_1(\kappa^{-1})$.

Composing  $J$ with $d$ yields a new real operator
\begin{equation}
d^{c}:=i(\overline{\partial}-\partial ).\label{n5}
\end{equation}
The formal adjoint of $d^{c}$ is denoted by $\delta^{c}$ and is given by
\begin{equation}
\delta^{c}=-\star d^{c}\star.\label{n6}
\end{equation}
On a K\"{a}hler manifold $d,d^{c}$ and $\delta,\delta^{c}$ are related
by the so-called K\"{a}hler identities
\begin{equation}
\hphantom{-}\delta=[\Lambda,d^{c}],~~~~-\delta^{c}=[\Lambda,d]\label{n7}
\end{equation}
\begin{equation}
-d=[L,\delta^{c}],~~~~\hphantom{-}d^{c}=[L,\delta]\label{n8}
\end{equation}
where $L$ is the algebraic operation
\begin{equation}
L\varphi=\omega\wedge\varphi\label{n9}
\end{equation}
and $\Lambda$ is its adjoint (contraction with $\omega$).  We note
\begin{equation}
[\Lambda,L]=(m-r){\bf 1}\label{n10}
\end{equation}
on $r$-forms.  Finally,
the Laplace-Beltrami operator $\Delta=\delta d$ on functions
may  be re-expressed in the useful form
\begin{equation}
\Delta f=-\Lambda dd^{c}f=-(\omega,dd^{c}f)~ .\label{n11}
\end{equation}

If ${\cal V}\to M$ is a holomorphic vector bundle over a complex manifold
$M$, its sheaf of sections will be denoted by $\O ({\cal V})$.
We define the {\em projectivization}  of ${\cal V}$
by $\bp ({\cal V})=({\cal V}-{\bf 0})/{\bc^{\times }}$, where
${\bf 0}$ is the zero section; notice
 that this differs from  a competing
convention which  replaces ${\cal V}$ with its  dual
${\cal V}\*$  on the right-hand side.
Finally, depending on the context, we will
use $\O (k)$ to
 denote either the degree $k$ holomorphic line bundle on $\bcp_m$,
or its sheaf of holomorphic  sections.

\subsection{Anti-self-duality}\label{asd}
On an oriented Riemannian 4-manifold $(M,g)$,
the bundle ${\wedge^{2}}$ of 2-forms
breaks up as
\begin{equation}
{\wedge^{2}}={\wedge^{+}}\oplus{\wedge^{-}}~,\label{n12}
\end{equation}
where $\wedge^{\pm}$ is the eigenspace of the Hodge operator $\star$
with eigenvalue $\pm 1$. We will call
 ${\wedge^{+}}$  the bundle of self-dual (SD) 2-forms and ${\wedge^{-}}$
the bundle of anti-self-dual (ASD) 2-forms. This decomposition
is conformally invariant, in the sense that it
is invariant under conformal rescalings $g\to e^u g$.

The decomposition (\ref{n12}) allows us to define  differential
operators $ d^{\pm}: {\cal E}^1\to {\cal E}(\wedge^{\pm})$ by following
the exterior derivative with projection $\wedge^2\to\wedge^{\pm}$.
Since a closed anti-self-dual form is automatically harmonic,
the following useful
vanishing result is an immediate consequence of Hodge theory:
\begin{propn}
If $M^4$ is compact and $\beta \in {\cal E}^1(M)$,
$d^{+}\beta =0 \Leftrightarrow d\beta =0$. {\hfill \rule{.5em}{1em} \\}
\label{sd2}
\end{propn}

Applying (\ref{n12}) to the curvature
operator ${R^{ab}}_{cd}: \wedge^2 \to \wedge^2$
results in a block-matrix decomposition

\begin{center}
\begin{tabular}{cccc}
\hphantom{~$\wedge_-$}& ~$\wedge^{+\ast}$~~ & ~~~~~$\wedge^{-\ast}$&
\hphantom{$\wedge_-$}\\
\end{tabular}

\begin{tabular}{c|c|c|}
\cline{2-3}&&\\
$\wedge^+$&$W_++\frac{s}{12}$&$\Phi$\\ &&\\
\cline{2-3}&&\\
$\wedge^-$&$\Phi$ & $W_-+\frac{s}{12}$\\&&\\
\cline{2-3}
\end{tabular}
\hphantom{$\wedge^-$}
\end{center}
where  $W_{\pm}$ are trace-free and
 $2\Phi$ is the trace-free part of the Ricci curvature $r$.
If $W_+=0$, the metric $g$ is said to be {\em anti-self-dual}, or ASD.
Since the Weyl curvature $W=W_++W_-$ is precisely the conformally invariant
piece of the curvature tensor $R$,  the  anti-self-duality condition
is invariant under conformal rescalings $g\to e^u g$; thus
it  makes sense to speak of ASD conformal (classes of) metrics.

 If the Riemannian manifold $(M,g)$ is actually
K\"{a}hler, and is given its canonical orientation,
 the decompositions (\ref{decomp})
and (\ref{n12}) are compatible in the sense that
\begin{equation}
{\wedge^{+}_{{\Bbb C}}}={\Bbb C}\omega\oplus \wedge^{0,2}\oplus
{\wedge^{2,0}}\label{n13}
\end{equation}
and
\begin{equation}
{\wedge^{-}_{{\Bbb C}}}={\wedge^{1,1}_{0}},\label{n14}
\end{equation}
where
${\wedge^{1,1}_{0}}=\{ \varphi\in {\wedge^{1,1}}~|~\omega\wedge\varphi =0\}$
is the bundle of ``primitive'' (1,1)-forms.
But, as a  consequence of (\ref{n3}),  the curvature operator of a
K\"ahler manifold is in $\mbox{End} (\wedge^{1,1})$; thus, the upper left-hand
block of the curvature operator must just
be a multiple of $\om\otimes \om^{\sharp}$.
This  immediately leads to the following:

\begin{propn} {\rm \cite{fla}\cite{G}}
In complex dimension 2, a
 K\"{a}hler metric $g$ is anti-self-dual iff it is
scalar-flat
 $(s \equiv 0)$.  {\hfill \rule{.5em}{1em} \\}    \label{sd1}
\end{propn}

We conclude this section with a  closely related observation.
If $\varphi$ is any form of type (1,1),  we can write
\begin{equation}
\varphi={\textstyle \frac{1}{2}}~(\tr\varphi)\omega+\varphi_{0}\label{n15}
\end{equation}
where
\begin{equation}
\tr\varphi :=(\varphi,\omega)\label{n16}
\end{equation}
and $\varphi_{0}\in{\wedge^{1,1}_{0}}$.  Applied to the Ricci form, this yields
\begin{equation}
\rho={\textstyle \frac{1}{4}}~s\omega+\rho_{0},\label{n17}
\end{equation}
so that
\begin{equation}
\rho~\mbox{is ASD}~~\Longleftrightarrow~~s=0~ .\label{n18}
\end{equation}
In particular, the Riemannian connections on
 $\kappa =\wedge^{2,0}$ and $\wedge^+ \cong\kappa \oplus {\Bbb R}$
are ASD iff the K\"ahler manifold $(M,g)$ is scalar-flat.

\subsection{Admissible K\"ahler Classes}
The fact that the Ricci form $\rho$ of a K\"ahler manifold represents
$2\pi c_1$  in deRham cohomology  leads to serious  constraints on
the scalar curvature of  K\"ahler metrics. Indeed,
on a compact complex surface with K\"ahler form $\om$,
the integral of the scalar curvature $s$, henceforth called the
{\em total scalar curvature}, must be given by
\be
\int_Ms\vol= 4\pi c_1\cdot [\omega ] \label{tsc}\ee
as an immediate consequence of (\ref{n17}). Since the volume of $M$
is just $[\omega ]^{2}/2$, we conclude that the
{\sl average value of the scalar curvature is determined by the K\"ahler
class alone}. (In complex dimension $m$, the total scalar curvature
is similarly given by
$\int s\vol = 4\pi  c_1
 \cdot [\omega]^{m-1}/{(m-1)!}$,
 while the volume is $[\omega]^{ m}/m!$; thus the average scalar curvature
is still completely determined by the K\"ahler class.)
 In trying  to classify scalar-flat
K\"ahler surfaces, the first logical step is therefore to limit
ourselves to those K\"ahler classes with total scalar curvature zero.
This motivates the following definition:

\begin{defn} Let $M$ be a compact complex surface. A K\"ahler class
$[\omega ]\in H^{1,1}(M, {\Bbb R})$ will be said to be {\em admissible} iff
$c_1\cup [\omega ]=0$, i.e. iff the total scalar curvature
$\int s\vol$ vanishes for  K\"ahler  metrics in $[\omega ]$.
The set of admissible K\"ahler classes will be denoted by ${\cal A}_M\subset
H^{1,1}(M)$.
\end{defn}
If $c_1^{\br}=0$, any K\"ahler class is admissible, and
Yau's solution \cite{yau2} of the Calabi conjecture asserts that every
such class is represented by a unique Ricci-flat metric; moreover,
this Ricci-flat metric is the only scalar-flat metric in the class,
since (\ref{n18}) tells us that $\rho/2\pi$ is the unique
harmonic representative
of $c_1^{\br}=0$ when $s=0$. If, on the other hand, $c_1^{\br}\neq 0$,
 the set ${\cal A}_M$ of admissible K\"ahler classes is evidently an
open set in a hyperplane in $H^{1,1}(M, {\Bbb R})$.
However, this open set is often {\em empty}, as is shown by the following
arguments of  Yau \cite{yau1}; {\em cf.}
 \cite{G2} \cite{h1}.

\begin{propn} Let $(M, J, \omega)$
 be a compact K\"ahler m-manifold, and suppose that
$L\to M$ is a holomorphic line bundle such that $c_1(L) \cup [\omega]^{m-1}=0$.
Then either  $\Gamma (M, {\cal O}(L^{\ell }))=0$
$\forall {\ell }\neq 0$, or else $c_1^{\Bbb R}(L)=0$.\label{ox}
\end{propn}
\begin{proof}
Suppose that
$u\in\Gamma (M, {\cal O}(L^{\ell }))$, $u\not\equiv 0$.
Using Poincar\'e duality,
the volume of the zero locus $u=0$,
counted with multiplicity, must equal $(c_1(L^{\ell }) \cup
[\omega]^{m-1})/(m-1)!=0$.
Thus $u\neq 0$, and $L^{\ell }$ is trivial.
\end{proof}

\begin{cor} Let $(M, J, \omega)$ be a compact K\"ahler manifold
of total scalar curvature zero. Then  either $c_1^{\Bbb R}=0$, or else
$\Gamma (M, {\cal O}(\kappa^{\ell }))=0$ for all $\ell \neq 0$.
In particular, the   Kodaira dimension of $M$  is either $0$ or $-\infty$.
\label{yup}
 \end{cor}
\begin{proof} Since the hypothesis may be interpreted as
stating  that $[\om ]^{m-1}\cup c_1(\kappa)=0$, we may
apply  Proposition \ref{ox} with $L=\kappa$.
\end{proof}

\begin{thm} {\rm \cite{yau1}}
Let $(M, J)$ be a compact complex surface
which carries an admissible K\"ahler class $[\omega ]$. Suppose, moreover, that
$M$ is not covered by a complex torus or a K3 surface. Then
$M$ is a ruled surface. In particular, $H^{0,2}(M)=H^{2,0}(M)=0$,
 and $b^{+}(M)=1$. \label{sd3}
\end{thm}
\begin{proof} From the conclusions of Corollary \ref{yup},
the Kodaira-Enriques classification \cite{bpv}
allows us to conclude, in the first instance,  that $M$ is covered by a
complex torus or a K3 surface, or that, in the second instance,
 $M$ is either a ruled surface or ${\Bbb CP}_2$.
Since we also have $H^0 (M, {\cal O}(\kappa^{\ell }))=0$ for ${\ell }<0$,
it follows
that $M\not\cong {\Bbb CP}_2$, so that, in the second instance,
 $M$ is  a ruled surface--- i.e. $M$ is obtained from
the  projectivization ${\Bbb P}(E)\to \Sigma_{\bf g}$ of a
rank 2 holomorphic vector bundle $E$ over
a compact complex curve $\Sigma_{\bf g}$ by blowing up $|\tau (M)|$ points.
\end{proof}

While the biregular classification of ruled surfaces over a curve
$\Sigma_{\bf g}$ involves the classification
of holomorphic vector bundles on the curve, the {\em bimeromorphic}
classification
is extremely simple. In fact \cite{bpv},
{\sl every ruled surface over $\Sigma_{\bf g}$
can be obtained from $\Sigma_{\bf g}\times \bcp_1$ by blowing up and blowing
down.}
   To see this, first
observe that,  since
$E\otimes {\cal L}$ is generated by its sections when
${\cal L}$ is sufficiently positive,  there
are sections of ${\Bbb P}(E)\to \Sigma_{\bf g}$ passing through
any given point. Choose three distinct such sections, and  successively
blow up each of their points of intersection while at each step blowing down
the
proper transform of the fiber through it. In finitely many steps this leads us
to
a minimal model equipped with a projection to $\Sigma_{\bf g}$
which admits three disjoint sections, and it is then easy to see
that this model must be the product surface $\Sigma_{\bf g}\times \bcp_1$.

Although  this knowledge will not prove essential for our
purposes, we conclude this section by pointing out  that the
set of admissible K\"ahler classes can be described in extremely
concrete terms:

\begin{propn} Let $(M, \om_0)$ be a compact K\"ahler surface.
Then ${\cal A}_M$ is precisely the set of  real $(1,1)$-classes
$[\om ]\in H^{1,1}(M, \br)$
satisfying the following cohomological conditions:
\begin{description}
\item{\rm (i)} $[\om ]\cdot c_1 =0$;
\item{\rm (ii)} $[\om ]\cdot [\om_0] > 0$;
\item{\rm (iii)} $[\om ]^2 > 0$; and
\item{\rm (iv)} $[\om ]\cdot C > 0$ for every curve $C\subset M$.
\end{description}\label{nakai}
\end{propn}
\begin{proof}
If $[\om ]$ satisfies {(i)} and {(iv)}, the proof of Proposition \ref{ox}
shows that either $H^{2,0}=0$ or else $\kappa$ is trivial.

If the former happens, $H^{1,1}=H^2$, and hence $H^2(M, {\Bbb Q})$ is
dense in $H^{1,1}(M, \br)$; in particular, $M$ is algebraic.
 When $[\om ]\in H^2(M, {\Bbb Q})$,
$k[\om ]\in  H^2(M, {\Bbb Z})$ for some $k\in {\Bbb N}$, and,  by
Nakai's criterion \cite{bpv},
$k[\om ]$ is thus a K\"ahler class iff {(iii)} and {(iv)} hold.
Thus the  cone of K\"ahler classes and the cone determined by conditions
{(iii)} and {(iv)} intersect $H^2(M, {\Bbb Q})$ in the same set.
But  both these cones are open and convex; since they contain the
same set of rational points, and since the open boxes with rational corners
form a basis for the topology of the Euclidean space $H^{1,1}(M, \br)=
H^2(M, \br )$,  they  must therefore
coincide. Conditions {(iii)} and {(iv)}  are therefore equivalent to
the  class $[\om ]$ being K\"ahler, and in particular imply {(ii)}. Condition
{(i)} is thus the only additional condition needed to assure a
class is admissible.

If, on the other hand, $\kappa$ is trivial, $M$ is either a torus or
a K3 surface. In the torus case every class is uniquely represented by
a form with constant coefficients, so that {(ii)} and {(iii)}
are easily seen to be necessary and sufficient for a class to be
K\"ahler (and automatically admissible). For the K3 case, the claim follows
from
   Todorov's surjectivity of the refined period map \cite{bpv}.
\end{proof}

\subsection{Holomorphic Vector Fields and Scalar Curvature}
\label{key}

The  space
$\Gamma (M, \O (T^{1,0}M))$ of
holomorphic vector fields on a complex manifold $(M,J)$
is equipped by the
Lie bracket with the structure of  a complex
Lie algebra,  denoted by $a(M)$. If $M$ is compact,
this is precisely the Lie algebra
of the group of biholomorphisms of $(M,J)$, since
a  vector
field $\Xi$ of type (1,0)  is holomorphic
iff  ${\pounds}_{\Re \Xi}J=0$,
where ${\pounds}$ denotes the Lie derivative.
In particular, if  $g$ is a K\"ahler metric
on $M$,  the Lie algebra $\imath(M,g)$ of real Killing fields
is canonically identified with a real sub-algebra of
$a(M)$ because the K\"ahler form,
being the unique harmonic representative of its deRham class, is
automatically invariant under the isometry group of $M$.
If, in addition,
 the scalar curvature of $M$ is constant, the following remarkable
result of  Lichnerowicz   \cite{lich},
which generalizes work of Matsushima \cite{mat},
says that, modulo parallel fields,  $a(M)$ is in fact just
the complexification
of  $\imath(M,g)$:

\begin{propn} {\rm (Matsushima-Lichnerowicz Theorem)}
If $(M, J, g)$ is a compact K\"ahler manifold of
constant scalar curvature,
$a(M)$ is the direct sum of the space of parallel
(1,0)-vector fields and the space of
vector fields of the form $(\bar{\partial}f)^{\sharp}$ where $f$ is
any (complex) solution
of the equation
\be\Delta^{2}f+2(dd^{c}f, \rho)=0.\label{vf2}\ee
Moreover,  a solution $f$ of (\ref{vf2}) corresponds to a Killing field iff
it is purely imaginary ($\Re f=0$). In particular, the algebra $a(M)$
is reductive (semi-simple plus Abelian), and the identity component of
the  group of biholomorphisms of $M$ has  a compact real form.
 \label{vf1} {\hfill \rule{.5em}{1em} \\}
\end{propn}

In the above proposition, $\sharp$ denotes, as always, the inverse of
$$\flat :{\Bbb C}\otimes TM\to {\Bbb C}\otimes
T^{\ast}M : X\mapsto g(X, \cdot )~,$$ and  in particular induces a
complex-linear isomorphism $T^{0,1\ast}M\stackrel{\sim }{\to} T^{1,0}M$.

\begin{defn} Let $M$ be a compact complex surface.
We will say that the {\em Matsushima-Lichnerowicz obstruction vanishes}
if the identity component of
the  group of biholomorphisms of $M$ has  a compact real form. \label{matlic}
\end{defn}

While the  Matsushima-Lichnerowicz Theorem gives us an important obstruction
to the existence of constant-scalar curvature K\"ahler metrics on a
compact complex manifold $M$ in terms of the algebra
$a(M)$ of holomorphic vector fields, a more subtle
such  obstruction  was later
discovered by Futaki \cite{fut}\cite{foot}.
The {\em Futaki character}
 ${\cal F}(\cdot , [\omega ]): a(M)\rightarrow {\Bbb C}$
 is defined by
\begin{equation}
{\cal F}(\Xi ,  [\omega ])= \int_{M}\Xi (\phi_{\om })\vol
\label{vf3}
\end{equation}
where $\Xi $ is any holomorphic vector field and $\phi_{\om }$ is the Ricci
 potential:
\begin{equation}
\rho=\rho_{ \rm H}+dd^{c}\phi_{\om }\label{vf4}
\end{equation}
where $\rho_{ \rm H}$ is harmonic and $\phi_{\om }$ is $C^{\infty}$, real
and normalized
so that
\begin{equation}
\int \phi_{\om }\vol=0.\label{vf5}
\end{equation}
Notice, by taking the trace of (\ref{vf4}), that
\begin{equation}
s=\mbox{constant}~\Longleftrightarrow~\phi_{\om }=0.\label{vf6}
\end{equation}

The most remarkable property of the Futaki character
${\cal F}(\cdot , [\omega ])$, implicit in our
notation but not evident from the definition,  is
\cite{cal2} \cite{fut} that it depends
 only upon the K\"{a}hler
{\em class}; for this reason it is sometimes referred to as the Futaki
invariant.  From (\ref{vf6}), it now follows immediately   that the  vanishing
of ${\cal F}(\cdot , [\omega ])$ is a necessary
condition for the  K\"{a}hler class $[\omega ]$
to contain a representative with constant
scalar curvature.

We now give a  way of rewriting the
Futaki invariant that is particularly useful when  $\Xi =
2(\bar{\partial}f)^{\sharp}$
for a complex-valued function $f$ on $M$ with $\int f\vol =0$.
Our calculations will actually work, however, for an
arbitrary holomorphic vector field, provided we define
its {\em holomorphy potential} $f$ by $f:= \b{\partial }\* {\bf G}\Xi^{\flat}=
{\bf G}
\bar{\partial }\* \Xi^{\flat}$, where ${\bf G}$ is the Green's operator of the
Hodge Laplacian. Now notice that the Ricci potential can be written
in  terms of the Green's operator and the scalar
curvature as $\phi_{\om }= -\frac{1}{2}{\bf G}s$. Hence
\bea {\cal F}(\Xi ,  [\omega ])&=&\int_{M}\Xi (\phi_{\om })\vol =
\langle \Xi^{\flat},  d\phi_{\om }\rangle \\&=&
\langle \Xi^{\flat},  \bar{\partial }\phi_{\om }\rangle =
\langle \Delta {\bf G}\Xi^{\flat},  \bar{\partial }\phi_{\om }\rangle\\&=&
\langle  2{\bar{\partial }} { \bar{\partial }\*}
{\bf G}\Xi^{\flat},  \bar{\partial }\phi_{\om }\rangle =
\langle  2\bar{\partial }f,  \bar{\partial }\phi_{\om }\rangle\\&=&
\langle  f,  2\bar{\partial }\* \bar{\partial }\phi_{\om }\rangle =
\langle  f,  \Delta\phi_{\om }\rangle\\&=&
\langle  f,  -{\textstyle \frac{1}{2}}(s-s_H)\rangle =
  -{\textstyle \frac{1}{2}}\int_{M}f(s-s_{ \rm H})\vol
\\&=& -{\textstyle \frac{1}{2}}\int_{M}fs\vol \eea
where $s_{ \rm H}$ is the average value of $s$ on $(M,g)$.
Notice, incidentally, that the conclusion is insensitive to the
normalization $\int f\vol =0$  if the total scalar curvature $\int s\vol
=s_{ \rm H}\int \volm$ happens to vanish.

As a consequence we deduce  an innocuous-looking fact (cf. \cite{besse},
 Proposition 2.159) which will later turn out to be  surprisingly
 important:

\begin{propn} Let $(M, \om)$ be a compact K\"ahler manifold of
{\em constant} scalar curvature $s=c$. Then, for any closed
(1,1)-form $\alpha$ one has
$$\left.\frac{d}{dt}{\cal F}(\Xi ,  [\omega  +t\alpha])\right|_{t=0}=
\langle f\rho  , \alpha_H \rangle~,$$
where $f$ is the holomorphy potential of $\Xi$  as defined above,
and $\alpha_H$ is the harmonic part of $\alpha$.\label{besser}
\end{propn}
\begin{proof}
Since the $\cal F$ only depends on the
K\"ahler class, we might as well assume that $\alpha$ is harmonic.
Since
${\cal F}(\Xi ,[  \om + t\om])= (1+t)^m
{\cal F}(\Xi , [\om ])=0$ by the assumption that
$\om$ has constant scalar curvature, whereas the corresponding
right-hand side $\langle f\rho  , \om \rangle
=\frac{1}{2}\int f s\vol =-{\cal F}(\Xi , [\om ])$
 vanishes for the same reason,
 we may  assume that the harmonic form  $\alpha$ is
{\em primitive}. This assumption
has the effect that the volume form  of $\om (t)$ is
$$\vol (t)=\frac{(\om +t\alpha )^m}{m!}= \vol + \frac{t\om^{m-1}\wedge
 \alpha }{(m-1)!}
  + O(t^2)=\vol + O(t^2)~,$$ so that the Ricci form, being determined by the
volume form and $J$, similarly satisfies
$$\rho (t) =\rho + O(t^2)~ .$$  The normalization of the
holomorphy potential reads
 $\int f(t)\vol =O(t^2)$ for the same reason. Hence
\bea \left.\frac{d}{dt}{\cal F}(\Xi ,  [\omega  +t\alpha])\right|_{t=0}
&=&{ -\frac{1}{2}}
\left.\frac{d}{dt} \left[\int_{M}f(t)s(t)\vol (t)\right]\right|_{t=0} \\&=&
{\textstyle -\frac{1}{2}}\left.\int_{M}\frac{df}{dt}\right|_{t=0}c\vol -
{\textstyle \frac{1}{2}}\left.\int_{M}f\frac{d}{dt}\left[s(t)\vol (t))
\right]\right|_{t=0}
\\&=&{\textstyle -\frac{1}{2}}
\left.\int_{M}f\frac{d}{dt}\left[s(t)\vol (t)\right]\right|_{t=0}
\\&=&-{\textstyle\frac{1}{(m-1)!}}
\int_{M}f\left.\frac{d}{dt}\left[\rho (t) \wedge \omega^{m-1} (t) \right]
\right|_{t=0}
\\&=&-{\textstyle\frac{1}{(m-1)!}}
\int_{M}f\rho\wedge \left. \frac{d}{dt} \omega^{m-1} (t) \right|_{t=0}
\\&=&-{\textstyle\frac{1}{(m-2)!}}
\int_{M}f\rho  \wedge \omega^{m-2} \wedge \alpha \\&=&
\int_{M}f\rho  \wedge \star  \alpha
=\langle f\rho  , \alpha_H \rangle~ .
\eea
\end{proof}

\subsection{Twistor Spaces}
\label{twistors}

The Penrose correspondence \cite{AHS}\cite{besse}\cite{P}  is a dictionary
between
  anti-self-dual conformal
Riemannian 4-manifolds and a special class of complex 3-folds.
We will begin by briefly explaining how the translation works in each
direction.

If $(M,[g])$ is an  anti-self-dual
Riemannian 4-manifold, its {\em twistor space} is a complex 3-manifold $Z$
whose underlying  smooth 6-manifold is the total space of the
sphere bundle of  the rank-three  real vector-bundle of  self-dual 2-forms:
\bea S^2\to &Z&:= \{ \om\in{\wedge}_+~|~\|\om\| =\sqrt 2 \}\\
&\hphantom{\wp}\downarrow\wp& \\
&M&\eea
 We now give
$Z$ an almost-complex structure $J: TZ \to TZ$, $J^2=-1$, by
first observing that,  for each $x \in M$,
there is a natural one-to-one correspondence
 between $\wp^{-1}(x)$
and the set of $g$-orthogonal complex structures $\jmath: T_xM\to T_xM$
inducing
the  given orientation on $T_xM$; namely, in the spirit of (\ref{n4}),
 any such $\jmath$ corresponds to the 2-form $\omega_{\jmath}$ defined by
$$ \omega_{\jmath} (\xi , \eta) =  g ({\jmath}\xi, \eta)~ .$$
Since the Levi-Civit\`a connection
of $g$ induces a splitting $TZ=H\oplus V$ of the  tangent bundle of
$Z$ into horizontal and vertical parts, and we have a canonical
isomorphism $\wp_{\ast}: H\to  \wp^{\ast}TM$, we may define
$J_H: H\to H$, $J_H^2=-1$ by $J_H|_{\om_\jmath}:=\jmath$.
Since the fibers of $\wp$ are oriented metric 2-spheres, we
may also define $J_V:V\to V$, $J_V^2=-1$ to be rotation by
$-90^{\circ}$ in the tangent space of the fiber.
Defining $J=J_H\oplus J_V$ then makes $Z$ an almost-complex manifold.
 Quite remarkably, this
almost-complex structure is  conformally invariant.
Even more remarkably, it  is {\em integrable} because its
Nijenhuis tensor may be identified with $W_+$, and so vanishes
precisely  by the assumption that $(M,g)$ is anti-self-dual.
The fibers of $\wp$ have become ${\Bbb CP}_1$'s  with normal bundle
$\O (1)\oplus \O (1)$ in the complex manifold $Z$, while the fiber-wise
antipodal map $Z\to Z : \om\mapsto -\om$ has become a free anti-holomorphic
involution $\sigma : Z\stackrel{\bar{\O}}{\to} Z$.

Conversely, let $Z$ be a complex 3-fold with free anti-holomorphic involution
$\sigma : Z\stackrel{\bar{\O}}{\to} Z$, $\sigma^2=\mbox{id}_Z$,
and suppose that there is a smooth $\sigma$-invariant rational curve
in $Z$ with normal bundle $\O (1)\oplus \O (1)$. Let ${\Bbb C}M$ denote
the connected component of this curve in
the  space of all ${\Bbb CP}_1\subset Z$ with normal bundle bundle
$\O (1)\oplus \O (1)$. Invoking  \cite{K},  ${\Bbb C}M$ is a complex
4-manifold. Moreover, $\sigma$ induces an anti-holomorphic
involution  $\hat{\sigma}: {\Bbb C}M\to {\Bbb C}M$
for which  our original curve corresponds to a
fixed point; the fixed-point set of $\hat{\sigma}$ is therefore a (non-empty)
real-analytic 4-manifold, the obvious
 connected component of which we denote
by  $M$. There is then an anti-self-dual conformal class of metrics
$[g]$ on $M$ determined by requiring that a complex tangent vector
$\xi\in{\Bbb C}\otimes TM=T{\Bbb C}M|_M$ satisfies $g(\xi , \xi)=0$
iff  the  corresponding section of the normal bundle
of ${\Bbb CP}_1\subset Z$ has a zero.
If $Z$ actually arises by the construction of the preceding paragraph,
and if our initial curve is a fiber of $\wp$, this
exactly  reconstructs the given manifold $M$
and conformal structure $[g]$.

Now recall from \S \ref{asd} that
a  K\"ahler manifold of complex dimension  2 is anti-self-dual  iff
 its scalar-curvature vanishes.
(This might be considered   rather
remarkable  insofar as neither
the K\"ahler condition nor the scalar-curvature condition are themselves
conformally invariant, and yet their coincidence is reflected by a
conformally invariant property.)
The Penrose correspondence thus
 provides an unexpected link between  K\"ahler surfaces
  and complex 3-folds. As will
now be explained, the  speciality of the
metric  being scalar-flat K\"ahler is echoed by a specialty of
 the twistor space $Z$ in a manner simple enough
to allow us to study scalar-flat K\"ahler
geometry by using Kodaira-Spencer theory. The central result here is a
theorem of Pontecorvo:

\begin{propn} {\rm \cite{Pt}} Let $\wp : Z\to M$ be the twistor
fibration of a  (perhaps non-compact) anti-self-dual
conformal Riemannian 4-manifold $(M, [g])$. Suppose we are given
 a complex hypersurface $D_1\subset Z$
for  which the restriction   $\wp|_{D_1}$ of the twistor
projection is a diffeomorphism onto $M$,
and let
$J$ denote the complex structure on $M$  given by this
section of $Z$.  (Thus $\wp|_{D_1}$ becomes a biholomorphism of $D_1$ and
$(M,J)$.)
Let  $D_2=\b{D_1}$  denote the image of $D_1$ under the real structure
$\s : Z\to Z$, and let $D=D_1\cup D_2$. Then there exists a
metric $g$ in the conformal class $[g]$ such that $(M,g,J)$ is
K\"ahler iff the divisor line bundle  of $D$ is
isomorphic to the half-anti-canonical line bundle $\kappa^{-1/2}$
of $Z$. \label{wanna}
\end{propn}

 The implication $\Rightarrow$ is relatively straightforward;
indeed,  the  K\"ahler form $\om$, being parallel and self-dual, is a
 solution of the twistor equation, and so defines, via the Penrose transform,
 a holomorphic section of $\kappa^{-1/2}$ vanishing precisely at
$D$. For an analogous proof in the $\Leftarrow$ direction, cf. \cite{L3}.

This implies the following key result, originally
discovered in a somewhat different guise by
Boyer \cite{boyer}:

\begin{thm} \label{bythm}  Let   $\wp : Z\to M$ be the twistor
fibration of a  compact anti-self-dual
4-manifold $(M, [g])$, and suppose that $b_1(M)$ is even.  Let  $D_1\subset  Z$
be a complex hypersurface which meets every fiber in exactly one point.
Then, for any $c\in \br^+$,
 the conformal class $[g]$ contains a unique scalar-flat
K\"ahler metric of  volume $c$. Conversely, every scalar-flat
K\"ahler surface arises in this way.
\end{thm}
\begin{proof} Let $[D]$  denote the divisor line bundle of the hypersurface
$D=D_1\cup \s (D_1)$. Then
$c_1([D])=c_1(\kappa^{-1/2})$, so that $[D]\otimes \kappa^{1/2}$
is an element of $\mbox{Pic}_0(Z)$. By the Penrose transform,
and using Proposition \ref{sd2},
$H^1(Z, \O)=\{ \beta \in {\cal E}^1(M)~|~d^+\beta =0\}/d{\cal E}^0(M)
=H^1(M, {\Bbb C})$,
so that one has
$\mbox{Pic}_0(Z)=H^1(Z, {\Bbb C}^{\times })=H^1(M, {\Bbb C}^{\times })$;
in other words, every topologically trivial holomorphic
line bundle on $Z$ admits a compatible
flat connection, and these all come from the
base. In particular, $[D]\otimes \kappa^{1/2}$ admits a flat
${\Bbb C}^{\times}$-connection, and is the pull-back of
a flat ${\Bbb C}^{\times}$-bundle
on $M$. In particular, $[D]\otimes \kappa^{1/2}$ is
trivial on $\wp^{-1}(U)$ for any sufficiently small open set $U\subset M$.
By Proposition \ref{wanna},  $g$ is therefore
 conformal to  a K\"ahler metric on a sufficiently small open set $U$.

We have thus shown that there are locally-defined smooth
functions $u\in {\cal E}_U$ for which the (1,1)-form $\om$ associated to
$(g,J)$ satisfies $0=d(e^u\om )= e^u(du\wedge\om + d\om )$.
Since any two local choices of $u$ differ by an additive constant,
the 1-form $\beta=-du$ is {\em  globally} defined on $M$, and
$$d\om = \beta \wedge \om ~ .$$
Because  the Fr\"ohlicher spectral sequence of any
compact complex surface  degenerates \cite{bpv}, the hypothesis
that $b_1(M)\equiv 0\bmod 2$ implies that
$H^1_{d}(M, \bc )= H^0(M, \Omega^1)\oplus \overline{H^0(M, \Omega^1)}$.
(A less elementary but deeper explanation of this decomposition stems from
the fact \cite{siu} that a compact complex surface admits  K\"ahler
metrics iff $b_1$ is even.)
Thus
the  closed real 1-form $\beta$ can be written as
$$\beta = \Re \alpha + df$$
for some holomorphic 1-form $\alpha$ and some smooth function $f$. Introducing
the conformally rescaled metric
$\hat{g}:=e^{-f}g$, we now have $d\hat{\om }=
\Re \alpha \wedge \hat{\om }$. But then
$$ 0= \int_M d(\alpha \wedge \hat{\om})=
-{\textstyle \frac{1}{2}} \int_M \alpha \wedge \b{\alpha} \wedge \hat{\om}=
{\textstyle \frac{i}{2}}  \|\alpha\|^2_{L^2,\hat{g}}~,$$
so that $\alpha =0$, and $\hat{g}$ is  K\"ahler.
Since $\hat{g}$ is also ASD, it is automatically  scalar-flat by
Proposition \ref{sd1}.
\end{proof}

\subsection{Deformation Problems}
For a compact manifold $M$ which admits a scalar-flat K\"{a}hler metric $g$,
a number of moduli problems are now obviously of interest:

\begin{description}
\item{(a)}
the moduli of scalar-flat K\"{a}hler metrics in the given K\"{a}hler
class;
\item{(b)}
 the moduli of scalar-flat K\"{a}hler metrics
for the given complex  structure;
\item{(c)}
 the moduli of scalar-flat K\"{a}hler metrics, with the  complex structure
allowed to vary; and
\item{(d)}
 the moduli of ASD conformal structures on $M$.
\end{description}
One  might also be tempted to  add the following:
\begin{description}
\item{(b$'$)}  the moduli of ASD Hermitian conformal structures for a given
complex
structure;
\item{(c$'$)}  the moduli of ASD Hermitian conformal
structures for some complex
structure.
\end{description}
However, as we saw in proving Theorem \ref{bythm},
a result of Boyer \cite{boyer}
 states that, because $b_1(M)$ is even,
(b$'$) and (c$'$) are  respectively equivalent to
(b) and
(c), so nothing is to be gained  by considering these problems
separately.

Of these problems, (a) can be tackled quite easily within the standard
framework of K\"{a}hler geometry, but for (b)--(d) very valuable information
comes from the twistor description.  As we saw in the previous section,
 the twistor space $Z$
of a scalar-flat K\"{a}hler surface $M$ is a complex 3-manifold equipped with a
real structure $\sigma$ and a $\sigma$-invariant divisor
$D$.  The complex structure of $Z$
completely determines
the conformal structure of $M$ while the divisor $D$
specifies the given complex structure on $M$.  Thus the moduli
problems (b)--(d)
correspond to the following problems in terms of $(Z,D)$:

\begin{description}
\item{(b$^{*}$)}  moduli of complex structures on $Z$ with $D$ as a fixed
$\sigma$-invariant
divisor;
\item{(c$^{*}$)}  moduli of complex structures on $Z$ which admit a
$\sigma$-invariant
divisor with divisor line-bundle isomorphic to  $\kappa^{-1/2}$;
\item{(d$^{*}$)}  moduli of complex structures on $Z$, which admit a
compatible
real structure $\sigma$.
\end{description}
Note that this point of view imposes different  equivalence relations
on the metrics occurring in the
different problems; in problems (a)--(b), two metrics of the same total volume
will be considered
equivalent iff they are literally {\em equal}, whereas in problems (c)--(d)
two metrics will be equivalent if they are in the same orbit of
 the diffeomeorphism group cross conformal rescalings.

Certainly an advantage of the starred formulation over the original one is that
 one can
appeal to the machinery of Kodaira-Spencer deformation theory to get
local information about the moduli spaces in terms of certain sheaf cohomology
groups of the twistor spaces. But from our point of view the
most significant advantage of this description is that the
cohomology groups involved in these distinct
problems are related by exact sequences;  once problem (b) is thoroughly
understood,
problems (c) and (d) can also be solved with relatively little further effort.

\subsection{Deformation Theory} \label{ks}

Let $Z$ be a compact complex manifold.
For us, a  {\em deformation} of $Z$ will consist of
 the  following:  a  ``parameter''
manifold ${\cal T}$ with basepoint $o$; a smooth manifold  ${\cal Z}$; a
proper  submersion
 $\varpi:{\cal Z}\rightarrow {\cal T}$;
an integrable fiber-wise complex structure on ${\cal Z}$;
and an identification of the
central fiber $\pi^{-1}(o)$ with $Z$.  If ${\cal T}'$ is another
manifold,  with
basepoint $o'$, and $\varphi:{\cal T}'\rightarrow {\cal T}$ is
a basepoint-preserving
smooth map, there is an
induced deformation $\varphi^{*}({\cal Z})\rightarrow
{\cal T}'$.  The deformation $\varpi:{\cal Z}\rightarrow {\cal T}$
is called {\em complete} if any other deformation can be induced
from it by a smooth map $\varphi$,
{\em versal} if, in addition,
 the derivative of $\varphi$ at $o'$ is always uniquely determined,
and {\em universal} if, in addition, the inducing $\varphi$ is always unique.
When a universal deformation of
$Z$ exists, a
 neighborhood of $o$ in the parameter space $\cal T$
 gives a  model for the
moduli space of complex structures on $Z$ in a neighborhood of the
given  structure.

If $\varpi:{\cal Z}\rightarrow {\cal T}$
is any  deformation in the above sense,
the {\em Kodaira-Spencer map} at $o\in {\cal T}$ is an
$\br$-linear map  ${\bf ks}: T_o{\cal T}\to H^1(Z_o, \Theta )$
obtained in \v{C}ech cohomology
by differentiating the transition functions of a fiber-wise complex
coordinate atlas on ${\cal Z}$. The first basic result of Kodaira-Spencer
theory is that a deformation is complete (respectively, versal) if
${\bf ks}$ is surjective (respectively, bijective). Notice that, by virtue of
its
definition,
the Kodaira-Spencer map behaves functorially under pull-backs.

The main results of Kodaira-Spencer theory
\cite{KS} assert that  any versal deformation
may be made into a holomorphic map $\varpi :{\cal Z}\to {\cal T}$
between complex manifolds (in an essentially
unique manner), and, more importantly,
 give sufficient conditions \cite{KNS} for
the existence of a versal or universal
deformation of $Z$ in terms of
the sheaf cohomology groups $H^{\jmath}(Z,\Theta)$, where
$\Theta ={\cal O}(T^{1,0}Z)$ is  the sheaf of holomorphic vector fields on $Z$.
These results may be summed up as follows:

\noindent
\begin{thm}
Suppose $H^{2}(Z,\Theta)=0$.  Then a (holomorphic) versal
deformation exists, with parameter space ${\cal T}$ an open
 neighborhood of $o=0$
in $H^{1}(Z,\Theta)$.  This deformation is universal if $H^{0}(Z,\Theta)=0$.
{\hfill \rule{.5em}{1em} \\}  \label{ks1}
\end{thm}

Unfortunately, this will not  suffice for our purposes,  because
we will be primarily interested
 in  deformations of complex manifolds {\em with real structure}. By
a real structure on a compact complex manifold $Z$, we always  mean
an anti-holomorphic involution of $Z$--- i.e. an
anti-holomorphic map $\sigma : Z\stackrel{\overline{\cal O}}{\to }Z$
such that $\sigma^2=\mbox{id}_Z$.
We will further  assume that $\sigma$
{\em acts freely}--- i.e. without fixed points.
This in particular means that $Z/\sigma$
is a smooth manifold, and  while the complex structure tensor
$J$ of $Z$  cannot descend to the quotient, the unordered
pair $\{ J,-J\}$ {\em is} globally well-defined downstairs.
We therefore introduce the following
concept:

\begin{defn} A {\em semi-complex manifold} is a smooth manifold
$P$, together with a 1-dimensional sub-bundle $L\subset \mbox{End}(TP)$
such that, in a neighborhood of
any point $x\in P$, $L$ is spanned by an integrable
 complex structure  $J$.
\end{defn}

Of course, near any point there are then exactly 2 choices of the
complex structure spanning $L$--- if $J$ is one, $-J$ is the other.
If $P^{2m}$ is a semi-complex manifold, we can thus equip $P$ with an atlas for
 which all the transition functions are either holomorphic or anti-holomorphic
diffeomorphisms of domains in ${\Bbb C}^m$; and
conversely, any manifold equipped with such an atlas has an induced
semi-complex structure.

\begin{example} Let $P$ be a smooth, unoriented surface, and let
$[g]$ be a conformal class of Riemannian metrics on $P$. Then
$[g]$ determines a unique semi-complex structure on $P$.
\end{example}

\begin{example} Let $Z$ be the twistor space of a half-conformally-flat
Riemannian 4-manifold $(M, g)$, and let $Z$ be its twistor space.
Let
$\sigma : Z\stackrel{\overline{\cal O}}{\to }Z$ be its real structure,
 acting  on the fibers of  $\wp : Z\to M$ by the antipodal map.
Then $P:=Z/\sigma$ is a semi-complex manifold. Notice, incidentally,
 that $P\to M$
is an ${\Bbb RP}^2$-bundle.
\end{example}

The following observation will be as  crucial as it is
 trivial:
{\em every semi-complex manifold $P$ is double-covered by a
complex manifold in a manner that makes the non-trivial deck transformation
a free anti-holomorphic involution.}
Indeed, one simply takes the cover to consist of
the  elements  $J\in L\subset \mbox{End}(TP)$ such that $J^2=-1$.
Thus, our basic example  $P=Z/\sigma$ of a semi-complex manifold, where
$Z$ is complex and $\sigma : Z\to Z$ is a free anti-holomorphic involution,
actually represents the general case.

If $Z$ is a complex manifold, the sheaf of holomorphic vector fields
is, as mentioned above,
 denoted by $\Theta :=\O (T^{1,0}Z)$. However, we may identify
the underlying real vector bundle of $T^{1,0}Z$ with the real tangent bundle
$TZ$ by $2\Re : T^{1,0}Z\to TZ: \Xi\mapsto \Xi + \b{\Xi }$, and in the
process we identify $\Theta$ with the sheaf
$$\Re \Theta : =\{ \xi \in {\cal E}(TZ)~|~{\pounds}_{\xi}J=0\}$$
of ``real holomorphic'' vector fields; of course, this only identifies them
as  sheaves of  real Lie algebras. The interesting observation is
 that $\Re\Theta$ is {\em exactly the same} for the conjugate complex manifolds
$(Z, J)$ and $(Z, -J)$, and is  thus even  well defined
on a semi-complex manifold.
This, of course, happens
 precisely because $\Re\Theta$ is the sheaf of infinitesimal automorphisms
of the semi-complex structure.

If we repeat our previous definitions of deformations and versality for
semi-complex manifolds, with the fiber-wise structures only required
to be semi-complex instead of complex, we immediately get the following
result:

\begin{propn} Let $P$ be a compact semi-complex manifold such that
$H^2(P, \Re\Theta)=0$. Then there exists a versal deformation of $P$
with  a neighborhood of
$0\in H^1(P, \Re\Theta)$ as parameter space.
This deformation is universal if $H^0(P, \Re\Theta)=0$.
\end{propn}
\begin{proof} The  Forster-Knorr  power-series proof  \cite{fok} of Theorem
\ref{ks}
goes through without any essential  changes.
\end{proof}

\begin{lemma} Let $Z$ be a complex manifold with free anti-holomorphic
involution
$\sigma : Z\to Z$, and let $P=Z/\sigma$ be the associated semi-complex
manifold.
Then $H^j(Z, \Theta)= H^j(P, \Re\Theta)\otimes_{\br } \bc$.
\end{lemma}
\begin{proof} There are arbitrarily fine covers $\cal V$ of $Z$ which are
equivariant under $\sigma$, and any such cover
descends to a cover $\cal W$ of $P$.
 For any such cover $\cal V$, $\sigma$ acts as on
$\check{H}^j({\cal V}, \Theta)$ via the anti-liner map
$\{ f_{\alpha\cdots \beta}\}\mapsto
\{ \sigma\* \overline{f_{\alpha\cdots \beta}}\}$,
and the fixed-point set of this action can be identified with
$\check{H}^j({\cal W}, \Re\Theta)$. Hence $\check{H}^j({\cal V}, \Theta)=
\check{H}^j({\cal W}, \Re\Theta)\otimes_{\br } \bc$. The lemma now follows by
taking direct limits.
\end{proof}

\begin{thm} Let $Z$ be a compact complex manifold with $H^2(Z, \Theta)=0$.
Suppose that $\sigma : Z\to Z$ is an anti-holomorphic involution without
fixed points. Then $\sigma$ can be extended as an anti-holomorphic
involution $\sigma_{\cal Z}:{\cal Z}\to {\cal Z}$
of the total space  of the versal deformation of
$Z$ which covers an anti-holomorphic involution
 $\sigma_{\cal T}:{\cal T}\to {\cal T}$ of the base.
The fixed-point set of  $\sigma_{\cal T}$ is a totally real
subspace ${\cal T}_{\sigma}$
of real dimension $h^1(Z, \Theta )$, and the restriction of
$\varpi :  {\cal Z}\to {\cal T}$ to this subspace is a versal deformation
of $({\cal Z}, \sigma )$.
\end{thm}

\begin{remark} When $H^0(Z, \Theta )=0$, this is an immediate consequence
\cite{DF} of the
universal property of the versal deformation.
\end{remark}

The importance of real deformations stems from the following observation,
the essence of which was discovered by Penrose \cite{P}:

\begin{thm} Let $\varpi : {\cal Z}\to {\cal T}$ be a deformation
of the twistor space $Z=Z_o$ of a compact  ASD conformal Riemannian 4-manifold
$(M, [g])$. Suppose, moreover, that ${\cal Z}$ is equipped with a
fiber-wise anti-holomorphic involution which restricts to
the twistor real structure on $Z_o$. Then there is a neighborhood
${\cal U}$ of $o\in {\cal T}$  and a  family  of ASD  Riemannian metrics $g_t$
on $M$, depending smoothly on $t\in {\cal U}$, such that
$Z_t=\varpi^{-1}(t)$ is biholomorphic to the twistor space of
$(M, [g_t])$.\label{Pen}
\end{thm}
\begin{proof} The point is that the normal bundle $\nu$ of a twistor fiber
$C=\wp^{-1}(x)\in Z$ is isomorphic to the bundle $\O (1)\oplus \O (1)$
on $\bcp_1$, and so satisfies $H^1(C, \nu )=0$. By Kodaira's
stability theorem \cite{K}, the complete analytic family generated by the
twistor fibers is stable under deformations. Because $\O (1)\oplus \O (1)$
is a rigid bundle on a rigid manifold, we have a 4-complex parameter
family of $\bcp_1$'s with normal bundle $\O (1)\oplus \O (1)$
in any small deformation $Z_t$ of $Z$; the  $\sigma$-invariant curves
in these families then foliate a manifold containing
$Z_o$ and spread over an open neighborhood
 of $Z_o\subset {\cal Z}$. There is therefore a neighborhood of
$Z_o\subset{\cal Z}$ foliated by these curves and projecting
properly to a neighborhood ${\cal U}$ of $o\in {\cal T}$.
Invoking  the inverse twistor correspondence described in \S\ref{twistors}
finishes the proof.
\end{proof}

In short, the moduli space for problem (d) is locally the same as the
moduli space for semi-complex structures on $Z/\s$. In order to
attack the moduli problems (b) and (c), we shall instead
 require  two `relative
versions' of the above deformation theory.  Suppose  we  are given a
compact complex manifold $Z$ and
a nonsingular complex hypersurface $D$ in
$Z$.  A deformation of $(Z,D)$ is given by the following:  a deformation
$\varpi:{\cal Z}\rightarrow {\cal T}$ of $Z$ and a deformation ${\cal D}
\rightarrow
{\cal T}$ of $D$, together with  a commutative diagram
\begin{eqnarray}
{\cal D}                       & \longrightarrow & {\cal Z}   \nonumber \\
\displaystyle\downarrow&                 &
\displaystyle\downarrow{\varpi} \nonumber \\
{\cal T}                              & =               &
 {\cal T} \label{ks2}
\end{eqnarray}
which restricts to the inclusion of $D$ in $Z$ at the central fibers.

A deformation of $(Z,D)$ with {\em fixed divisor} $D$ is a deformation of
the above type with ${\cal D}\cong D\times {\cal T}$, and $\varpi |_{\cal D}$
corresponding to
projection on the
second factor.
The notions of versal and universal deformations exist also for these
relative deformations and there are analogues of Theorem \ref{ks1}.  To
state them, let $\Theta_{Z,D}$ be the sheaf of holomorphic vector fields on
$Z$ that are tangent to $D$ along $D$; and let ${\Theta}_{Z}\otimes {\cal I}_D$
be the
subsheaf of vector fields which vanish along $D$.

Finally, we can modify all the above definitions so as to replace
$Z$ with a semi-complex manifold $P$ and $D$ with a
semi-complex submanifold $M\subset P$ of real codimension 2.
The same reasoning as before then yields

\noindent
\begin{thm}
Let $(P,M)$ be a semi-complex manifold with nonsingular
 semi-complex hypersurface,
and let $(Z, D, \sigma )$ be the complex manifold with hypersurface
and real structure which covers it.

(i)  Suppose that $H^{2}(Z,\Theta_{Z,D})=0$.
Then there are  versal deformations of $(P,M)$ and $(Z,D)$.
 Moreover, the parameter space for the former deformation is
a real slice in  that of the latter,
which may be taken to be a neighborhood of $0\in H^1(Z,\Theta_{Z,D})$.
These deformations are both
universal if $H^{0}(Z,\Theta_{Z,D})=0$.

(ii) Suppose that $H^{2}(Z,{\Theta}_{Z}\otimes {\cal I}_D)=0$.  Then there are
versal deformations of $P$ with fixed divisor $M$, and
of $Z$ with fixed divisor $D$.
Moreover, the parameter space for the former deformation is
a real slice in that  of the latter,
which may be taken to be a neighborhood of
$0\in H^1(Z,{\Theta}_{Z}\otimes {\cal I}_D)$.
These
deformations are both
 universal if $H^{0}(Z,{\Theta}_{Z}\otimes {\cal I}_D)=0$.    \label{ks3}
\end{thm}

Notice that these deformation problems are all related in that there exist
exact sequences:
\begin{equation}
0\rightarrow\Theta_{Z,D}\rightarrow\Theta_{Z}\rightarrow N_{D}\rightarrow 0
\label{ks4}
\end{equation}
where $N_{D}$ is the normal bundle of $D$ in $Z$, and
\begin{equation}
0\rightarrow{\Theta}_{Z}\otimes {\cal I}_D\rightarrow
\Theta_{Z,D}\rightarrow\Theta_{D}
\rightarrow 0.\label{ks5}
\end{equation}
Moreover, the induced long-exact  sequences exactly intertwine
the  Kodaira-Spencer maps of the deformation theories involved.
In particular, given a deformation
$\varpi : ({\cal Z}, {\cal D})\to {\cal T}$ of
$({ Z}, { D})$, there is a  Kodaira-Spencer map ${\bf ks}_{Z,D}\in
\mbox{Hom} (T_o{\cal T},  H^1(Z, \Theta_{Z,D}))$,
gotten by differentiating  the transition functions
of a fiber-wise complex atlas which sends open sets of $Z_t$ to $\bc^m$
and open sets of $D_t$ to $\bc^{m-1}\subset \bc^m$; this is obviously
related to
the Kodaira-Spencer maps of $\varpi : {\cal Z}\to {\cal T}$
and  $\varpi |_{\cal D}:{\cal D} \to {\cal T}$ by composition with
the  natural homomorphisms
$H^1(Z, \Theta_{Z,D})\to H^1(Z, \Theta_{Z})$ and $H^1(Z, \Theta_{Z,D})
\to H^1(D, \Theta_{D})$. In particular, if $H^2 (Z, {\Theta}_{Z}\otimes {\cal
I}_D)
=0$, the resulting surjectivity of $H^1(Z, \Theta_{Z,D})
\to H^1(D, \Theta_{D})$ implies that if
$\varpi : ({\cal Z}, {\cal D})\to {\cal T}$
is assumed to be a versal deformation of $({ Z}, { D})$, the
induced deformation $\varpi |_{\cal D}:{\cal D} \to {\cal T}$
is complete. A similar argument in the semi-complex case will feature
prominently in our proof of the Main Theorem.

\pagebreak
\setcounter{equation}{0}
\section{Deformations of Scalar-flat K\"{a}hler
Surfaces}\label{flag}

The standard treatment of these problems is as follows.  Given a scalar-flat
K\"{a}hler metric $g$, K\"{a}hler form $\omega$, normalized so that the
total volume is 1, we identify the tangent space to the space of volume-1
K\"{a}hler forms as
\begin{equation}
K=\left\{\varphi\in{\wedge^{1,1}}(M):d\varphi=0~
\mbox{ and }~\int~(\tr\varphi)
\vol=0\right\}.\label{d1}
\end{equation}
The derivative in the direction $\varphi$ of the scalar curvature is (cf.
\cite{besse}, Lemma 2.158(iii))
\begin{equation}
s'(\varphi)=\Delta (\tr\varphi)-2(\rho,\varphi)\label{d2}
\end{equation}
(where $\rho$ is the Ricci form as in (\ref{n4})).

If it is required to preserve the K\"{a}hler class then $\varphi\in K$ is
taken to have the form
\begin{equation}
\varphi=-dd^{c}f\label{d3}
\end{equation}
(for some real $C^{\infty}$ function $f$) and (\ref{d2}) reduces to
\begin{equation}
s'(f)=\Delta^{2}f+2(dd^{c}f, \rho)\label{d4}
\end{equation}
This we recognize as Lichnerowicz's differential equation (\ref{vf2}).
Invoking Proposition \ref{vf1}, we can thus
immediately solve problem (a):

\begin{propn}
The tangent space to the moduli space of
scalar-flat K\"{a}hler metrics in a given K\"{a}hler class is precisely
the space $\imath(M)^{\perp}\subset a(M)$
 of holomorphic vector fields orthogonal to the space $\imath(M)$ of
Killing fields.  {\hfill \rule{.5em}{1em} \\}
\label{dt1}

\end{propn}

Let us now turn to the more general problem (b).  For this, we shall study
the equivalent problem of deforming $Z$ with fixed divisor $D$.  Referring
to Theorem \ref{ks3} we see that the first task
is to identify ${\Theta}_{Z}\otimes {\cal I}_D$ and its cohomology groups.
Since the ideal sheaf of $D$ is isomorphic to
${\cal O}(-2):=\kappa_Z^{1/2}$, we have
${\Theta}_{Z}\otimes {\cal I}_D=\Theta(-2)$.  To study the
cohomology groups we use
the Penrose transform to relate them to data on $M$.  A straightforward
application of the  techniques
of \cite{be}\cite{bs} or \cite{h2}
yields:

\noindent
\begin{propn} Let $(M,g)$ be any anti-self-dual manifold.
For each $j=0,\ldots,3$, the Penrose
transform identifies $H^{j}(Z,\Theta(-2))$ with the ${j}$-th
cohomology group of
the complex
\begin{equation}
0\rightarrow{\wedge^{-}}(M)~ \stackrel{S}{\rightarrow}~{\wedge^{+}}(M)
\rightarrow 0
\end{equation}
where
\begin{equation}
S(\alpha)=d^{+}\delta\alpha + \Phi\alpha 
\end{equation}
for $\alpha\in{\wedge^{-}}(M)$. Here $\Phi : \wedge^-\to \wedge^+$ denotes
one-half the trace-free Ricci curvature, acting by $\alpha_{ab}\mapsto
\Phi^c_{b}\alpha_{ac}-\Phi^c_{a}\alpha_{bc}$. \label{dt0}
 {\hfill \rule{.5em}{1em} \\}
\end{propn}

\noindent{\bf Remark}.  The  operator $S$  is
 conformally invariant, provided that
$\alpha$ is   conformal weighted as follows:
  $\alpha
\mapsto e^{u/2}\alpha$ when $g\mapsto e^{u}g$.

\begin{cor} Let $(M,g)$ be a scalar-flat K\"ahler surface.
The Penrose
transform then
identifies $H^{j}(Z,\Theta(-2))$ with the ${j}$-th cohomology group of
the complex
\begin{equation}
0\rightarrow{\wedge^{-}}(M)~ \stackrel{S}{\rightarrow}~{\wedge^{+}}(M)
\rightarrow 0\label{d5}
\end{equation}
where
\begin{equation}
S(\alpha)=d^{+}\delta\alpha-{\textstyle
\frac{1}{2}}~(\rho,\alpha)\omega\label{d6}
\end{equation}
for $\alpha\in{\wedge^{-}}(M)$.  {\hfill \rule{.5em}{1em} \\} \label{dt2}
\end{cor}

While the  twistor theory predicts that the
operator $S$ completely  governs problem (b), it is perhaps not obvious
why this is so. Let us therefore digress for a moment in order to observe
 the kernel of $S$ can indeed be identified with space of
$\varphi$ given by (\ref{d1}) and (\ref{d2}).

\begin{thm}
The map
$$
K\rightarrow{\wedge^{-}}(M)
$$
given by $\varphi\mapsto\varphi_{0}$ (see (\ref{n15})) induces
an isomorphism of $\ker(s')$ with $\ker(S)$.    \label{dt3}
\end{thm}

\noindent{\bf Proof}.  We begin by noticing that, as a consequence of
 (\ref{n8}) and (\ref{n13}), the equation
\begin{equation}
d^{+}\delta\alpha=\lambda\omega\label{d11}
\end{equation}
is equivalent to the two equations
\begin{equation}
{\Lambda}d\delta\alpha={-\Lambda}\delta d\alpha=2\lambda\label{d9}
\end{equation}
and
\begin{equation}
d^{+} \Lambda d\alpha=0.\label{d10}
\end{equation}

Now suppose that $\varphi$ is in $K$ and write
$$
\varphi={\textstyle \frac{1}{2}}~(\tr\varphi)\omega+\varphi_{0}.
$$
We have to show that equation (\ref{d2}) implies that
$\varphi_{0}$ is in the kernel of $S$.  But $d\varphi =0$ is
equivalent to
\begin{equation}
d\varphi_{0}=-{\textstyle \frac{1}{2}}Ld\,(\tr\varphi)\label{d12}
\end{equation}
and so to
\begin{equation}
{\Lambda}d\varphi_{0}=-{\textstyle \frac{1}{2}}~[{\Lambda},L]d\tr\varphi
=-{\textstyle \frac{1}{2}}
{}~d(\tr\varphi)\label{d14a}
\end{equation}
by (\ref{n10}) (remember $\Lambda$ of any 1-form is zero).  This implies that
$\varphi_{0}$ satisfies equation (\ref{d10}).  On the other hand, by applying
${\Lambda}\delta$ to (\ref{d12}) we get
$$
\begin{array}{lcll}
{\Lambda}\delta d\varphi_{0} & = & - &
{\textstyle \frac{1}{2}}
{\Lambda}\delta Ld(\tr\varphi) \\[+12pt]
 & = & - &
{\textstyle \frac{1}{2}}{\Lambda}[\delta,L]d\tr\varphi-
{\textstyle \frac{1}{2}}{\Lambda}L\delta d(\tr\varphi) \\[+12pt]
 & = & & {\textstyle \frac{1}{2}}~
{\Lambda}d^{c}d(\tr\varphi)-
\delta d(\tr\varphi) \\[+12pt]
 & = & & {\textstyle \frac{1}{2}}~
\Delta (\tr\varphi)-\Delta (\tr\varphi)=-
{\textstyle \frac{1}{2}}~\Delta (\tr\varphi)
\end{array}
$$
where we have used the K\"{a}hler identities (\ref{n8}), (\ref{n10}) and
(\ref{n11}).
Hence by equation (\ref{d2})
$$
{\Lambda}\delta d\varphi_{0}=-(\rho,\varphi)=-(\rho,\varphi_{0})
$$
because $\rho$ is ASD (cf. (\ref{n18})), and this is equation (\ref{d9}) with
$\lambda={\textstyle \frac{1}{2}}~(\rho,\varphi_{0})$ as required.

   To go in the other direction we suppose $\alpha\in\ker(S)$
so that it satisfies
(\ref{d9}) and (\ref{d10}) with $\lambda=
{\textstyle \frac{1}{2}}~(\rho,\alpha)$.
Let $u$ be
the unique solution of $\Delta u=2(\rho,\alpha)$ with $\int u\vol=0$, and put
$\varphi={\textstyle \frac{1}{2}}~u\omega+\alpha$.  Then $\varphi$
automatically satisfies
(\ref{d2}): all that remains is to check $d\varphi=0$.

By Proposition \ref{sd2} and equation (\ref{d10}), ${\Lambda}d\alpha$
is $d$-closed.  Accordingly its Hodge decomposition takes the form
\begin{equation}
{\Lambda}d\alpha=h+dv\label{d13}
\end{equation}
where $h$ is a harmonic 1-form and $v$ is a $C^{\infty}$ function which is
unique if we insist that $\int v\vol=0$.  We claim that $h=0$.  Indeed
$$
||h||^{2}=\langle h,{\Lambda}d\alpha+dv\rangle=\langle h,-\delta^{c}\alpha
\rangle=-\langle d^{c}h,\alpha\rangle=0
$$
where we've used the K\"{a}hler identity (\ref{n7}) and the basic fact that on
a K\"{a}hler manifold any $\Delta$-harmonic form is also $\Delta^{c}$-harmonic.
Now if we compare (\ref{d13}) with (\ref{d14a}) and
the definition of $u$ we see
that $d\varphi=0$ iff $v=-{\textstyle \frac{1}{2}}~u$.
   To see that this is the
case we use (\ref{d13}) to compute
$$
\Delta(-2v)  =  -2\delta{\Lambda}d\alpha=2{\Lambda}d\delta\alpha=2(\rho,
\alpha) = \Delta u
$$
as required.  In the above we have used the K\"{a}hler identities and
equation (\ref{d9}).  {\hfill \rule{.5em}{1em} \\}

\begin{remark} Aside from the twistor-theoretic argument,
the relevance of the operator  $S$ to
 problem (b) can best be seen  by first  restating
 the problem as  problem (b$'$).  One then observes that
$S$ is the linearization of the operator  which sends a
Hermitian metric $g$, with associated 2-form $\om_g$, to
$W_{+g}(\om_g )\in {\cal E}(\wedge^+)$. From Boyer's calculations
\cite{boyer}
one then
reads off the fact that the kernel of this non-linear operator
 is precisely the space
of ASD Hermitian conformal classes, and, since $b_1(M)$ is even,
 these are all  represented by unit-volume scalar-flat K\"ahler metrics.
\end{remark}

\noindent
\begin{propn}
The operator $S$ of (\ref{d6}) is
elliptic with index equal to $-\tau(M)$, where $\tau(M)$ is the signature
of $M$. \label{dt2a}
\end{propn}
\begin{proof} Both statements depend only on the top-order term
$d^{+}\delta$ of
$S$. Now the kernel of  $d^{+}\delta$ is
the space $H^-$ of ASD harmonic 2-forms and similarly  the kernel of its
adjoint
$d^-\delta$ is the space $H^+$
of SD harmonic 2-forms. The proof that the symbol is an isomorphism
$\wedge^-\rightarrow \wedge^+$ is left to the reader.
\end{proof}

We now complete our analysis of problem (b) by identifying the cokernel of
$S$.

\begin{propn} Suppose that $M$ is not Ricci-flat.
The cokernel of $S$ can then be identified
with the space of $C^{\infty}$ functions $f$ which satisfy the following
conditions:

\noindent$\!\!\!$\begin{tabular}{lp{5.9in}}
(i) & the Lichnerowicz equation (\ref{vf2}):
$~\Delta^{2}f=-2(dd^{c}f, \rho)$; \\
(ii) & the orthogonality conditions $\langle f\rho,\alpha\rangle=0$
for all ASD
harmonic 2-forms $\alpha$.\label{crux}
\end{tabular}

\noindent In particular, if $M$
supports no non-parallel holomorphic vector
fields, then coker$(S)=0$.     \label{dt4}
\end{propn}

\noindent{\bf Proof}.  By the Fredholm alternative for elliptic operators,
the cokernel of $S$ can be identified with the kernel of the adjoint
$S^{*}$.  Now for any $\psi\in{\Lambda^{+}}(M)$ and $\alpha\in{\wedge^{-}}(M)$,
$$
\begin{array}{lcl}
\langle S^{*}\psi,\alpha\rangle & = & \langle\psi,S\alpha\rangle=\langle\psi,
d^{+}\delta\alpha-{\textstyle \frac{1}{2}}~(\rho,\alpha)\omega\rangle \\
 & = & \langle d^{-}\delta\psi,\alpha\rangle-
{\textstyle \frac{1}{2}}~\langle(\psi,\omega)
\rho,\alpha\rangle
\end{array}
$$
so
\begin{equation}
S^{*}\psi=d^{-}\delta\psi-{\textstyle
\frac{1}{2}}~(\psi,\omega)\rho.\label{d14}
\end{equation}
   To analyze the equation $S^{*}\psi=0$, we shall invoke (\ref{n13}) to write
\begin{equation}
\psi=f\omega+\chi\label{d15}
\end{equation}
(where $f$ is a $C^{\infty}$ function and $\chi$ lies in ${\wedge^{2,0}}
\oplus{\wedge^{0,2}})$.  We shall also need to write the operator $d^{-}
\delta$ in terms of $d$ and $d^{c}$.  This is an exercise involving the
K\"{a}hler identities (\ref{n7}).  Indeed, as an operator ${\wedge^{+}}
\rightarrow{\wedge^{-}}$,
\begin{eqnarray}
d^{-}\delta & = & {\textstyle \frac{1}{2}}~(1-\star)d\delta \nonumber \\
                 & = & {\textstyle \frac{1}{2}}~d\delta-
{\textstyle \frac{1}{2}}~\delta d
\nonumber  \\
                 & = &
{\textstyle \frac{1}{2}}~d[{\Lambda},d^{c}]-{\textstyle \frac{1}{2}}~
[{\Lambda},d^{c}]d \nonumber \\
                 & = &
{\textstyle \frac{1}{2}}~(d{\Lambda}d^{c}+d^{c}{\Lambda}d)+
{\textstyle \frac{1}{2}}~
({\Lambda}dd^c-dd^{c}{\Lambda}). \label{d16}
\end{eqnarray}
Consider the second bracketed term in (\ref{d16}).
For reasons of bidegree, it
annihilates $\chi$ in (\ref{d15}).  On the other hand
\begin{eqnarray}
d^{-}\delta(f\omega) & = & d^{-}\delta Lf \nonumber \\
                     & = & -d^{-}[L,\delta]f \nonumber \\
                     & = & -{\textstyle \frac{1}{2}}
{}~(1-\star)dd^{c}f \nonumber \\
                     & = & -dd^{c}f-{\textstyle \frac{1}{2}}
{}~\omega\Delta f \label{d17}
\end{eqnarray}
where we have used the K\"{a}hler identities, the relation (\ref{n14}) to
identify the ASD part of $dd^{c}f$ with its projection perpendicular to
$\omega$, and (\ref{n11}) to relate this to the Laplacian.  Combining
(\ref{d15}), (\ref{d16}) and (\ref{d17}), we obtain
\begin{equation} \label{d18}
S^{*}(f\omega+\chi)= d^-\delta \chi
-dd^{c}f-{\textstyle \frac{1}{2}}~\omega\Delta f-\rho f
\end{equation}
where we can also write
\begin{equation}
d^{-}\delta\chi={\textstyle \frac{1}{2}}~(d{\Lambda}d^{c}+d^{c}{\Lambda}d)\chi.
\label{d19}
\end{equation}

Suppose (\ref{d18}) vanishes.  Applying $dd^{c}$ (and using (\ref{d19})) we
find that $f$ satisfies condition (i) of the Theorem:
$$
\begin{array}{lcl}
0 & = & -~{\textstyle \frac{1}{2}}
{}~\omega{\wedge}dd^{c}\Delta f-\rho{\wedge}dd^cf \\
 & = & -~{\textstyle \frac{1}{2}}
{}~(\omega,dd^{c}\Delta f)\vol+(dd^{c}f, \rho)\vol \\
 & = & \left({\textstyle \frac{1}{2}}~\Delta^{2}f+(dd^{c}f, \rho)\right)\vol
\end{array}
$$
(the change of sign in the term in $\rho$ arises
 because $\rho$ is ASD cf. (\ref{n18})).
The orthogonality conditions (ii) are just the conditions that the equation
(\ref{d18})
\begin{equation}
d^{-}\delta\chi=dd^{c}f+{\textstyle \frac{1}{2}}~\omega\Delta f+\rho f
\label{d20}
\end{equation}
be soluble for $\chi$.  By the Fredholm alternative this equation is
soluble iff the right-hand side is orthogonal to the kernel of the adjoint
operator $d^{+}\delta$.  But we have already identified this kernel in the
Proof of Proposition \ref{dt2a} with the space $H^{-}$ of ASD harmonic 2-forms.
Since the inner product of the first two terms on the right-hand side
 with any such form
is zero we get condition (ii) of the Theorem.  The proof is completed by
noting that if $\chi$ satisfying  (\ref{d20}) exists, it is unique.  This is
because the kernel of $d^{-}\delta$ is $H^{+}={\Bbb C} \omega$ by
 Theorem \ref{sd3}, and by definition $\chi$ is orthogonal to
$\omega$.  {\hfill \rule{.5em}{1em} \\}
\vspace{12pt}

\begin{defn} Let $M$ be a compact complex surface.
If $\Xi$ is any holomorphic vector field on $M$, the {\em
restricted
Futaki invariant} of $(M, \Xi )$ is  defined to be the map
\bea \hat{\cal F}_{\Xi}: {\cal A}_M & \to & {\Bbb C} \\  ~~
 [\omega ] & \mapsto & {\cal F}(\Xi ,  [\omega ]) ~ .\eea\label{two}
\end{defn}
Here ${\cal A}_M:=\{ [\om ]\in H^{1,1}~|~[\om ]>0, c_1\cup [\omega ]=0\}$
again denotes the set of admissible K\"ahler classes.

\begin{thm} Let $(M, \omega )$ be a compact scalar-flat K\"ahler surface,
and let $Z$ be its twistor space. Assume that $M$ is not Ricci-flat. Then the
cohomology groups $H^{2}(Z,{\Theta}_{Z}\otimes {\cal I}_D)$,
$H^{2}(Z, {\Theta}_{Z,D})$, and $H^{2}(Z,{\Theta}_{Z})$
are all equal, and can be identified
 with the space of  holomorphic vector
fields $\Xi $ on $M$ such that $d\hat{\cal F}_{\Xi}|_{[\omega ]}=0.$
\label{dt5}\end{thm}

\begin{proof} Let us first observe that $M$ cannot carry a non-zero
parallel vector field. If it did, $g$ would locally be a Riemannian
product of the flat metric  and some other
 K\"ahler metric on  ${\Bbb C}$;
and since $s=0$, the second factor would also have to be flat.
Thus $g$ would itself be flat,
contradicting  the assumption that $\rho \not\equiv 0$.
Since
$g$ has constant scalar curvature and
there are now no parallel vector fields on $M$, we may therefore, by
 Proposition \ref{vf1},
write each holomorphic vector field
$\Xi$ in the form $2(\bar{\partial} f)^{\sharp}$ for a unique $f$
satisfying the Lichnerowicz equation (\ref{vf2})
 and $\int f\vol =0$.

Now, in accordance with
Definition \ref{two} above, the   restricted
Futaki invariant  $\hat{\cal F}_{\Xi}$  is just the restriction of
${\cal F}(\Xi,\cdot )$,  defined by  (\ref{vf3}), to the admissible
K\"ahler classes ${\cal A}_M\subset H^{1,1}$.
The tangent space of ${\cal A}_M$ at $[\om ]$ is  just
the $\cup$-orthogonal complement of
$\rho$ in the harmonic  (1,1)-forms; but since the Futaki invariant
vanishes for all multiples of $[\om ]$, we might as well  restrict ourselves
to admissible classes of {\em fixed volume}, which corresponds to
  cutting the tangent space
down to the $L^2$-orthogonal complement of
$\rho$ in the closed ASD 2-forms $H^{-}$.
If $\Xi =2(\bar{\partial} f)^{\sharp}$, $\int f\vol =0$, then, by
Proposition \ref{besser} we have
$$\frac{d}{dt}{\cal F}(\Xi , [\om + t\alpha ])= \langle f\rho,\alpha\rangle .
$$
However, if $C$ is any constant and $\hat{f}:=f+C$, this becomes
$$\frac{d}{dt}{\cal F}(\Xi , [\om + t\alpha ])
=\langle\hat{f}\rho,\alpha\rangle-C\langle\rho,\alpha
\rangle
$$
so that the right-hand side  is independent of the
representative $\hat{f}$ provided  $\langle \rho,\alpha\rangle =0$.
Thus $$d\hat{\cal F}_{\Xi}(\alpha)=\langle f\rho,\alpha\rangle
{}~~\forall \alpha\in H^-~s.t.~\langle \rho,\alpha\rangle =0$$
for any $f$ with $\Xi =2(\bar{\partial} f)^{\sharp}$, independent
 of any statement concerning $\int f\vol$.
 On the other hand, if $\Xi =2(\bar{\partial} f)^{\sharp}$
is a holomorphic  vector field with $d\hat{\cal F}_{\Xi}|_{[\om ]}=0$,
there is exactly one $C$ for which $\hat{f}:=f+C$
satisfies
 $\langle\hat{f}\rho,\rho\rangle=0$, since
$\langle\rho,\rho\rangle=-4\pi^{2}c_{1}^{2}>0$.
This allows us to identify the space of holomorphic vector fields
$\Xi$ satisfying $d\hat{\cal F}_{\Xi}|_{[\om ]}=0$
with the space of solutions $f$ of the Lichnerowicz equation
such that $\langle f\rho,\alpha\rangle=0$ for all $\alpha \in H^-$.
Using Corollary \ref{dt2} and Theorem \ref{crux}, this in turn identifies
$H^{2}(Z,{\Theta}_{Z}\otimes {\cal I}_D)$ with the space of
holomorphic vector fields $\Xi$ on $M$ such that $d\hat{\cal F}_{\Xi}|_{[\om
]}=0$,
as promised.

Now, using Serre duality, we observe that, since $M$ is ruled,
  $H^2(M, \Theta_M)\cong H^0 (M, \Omega^1 ({\kappa}_M ))=0$ because
 $\Omega^1 (\kappa )$  becomes
 $ \O (-2 )\oplus \O (-4 )$ when restricted to a smooth rational curve
with trivial normal bundle. Since $[\omega ]$ has  total scalar curvature 0,
we also have   $H^2(M, \O (\kappa^{-1}_M)) \cong H^0 (M, \O (\kappa^2_M))=0$
by Corollary \ref{yup}. The isomorphisms
$$H^{2}(Z,{\Theta}_{Z}\otimes {\cal I}_D)\cong
H^{2}(Z, {\Theta}_{Z,D}) \cong H^{2}(Z,{\Theta}_{Z})$$
now follow immediately from the short exact sequences (\ref{ks4}) and
(\ref{ks5}),
since $D$ consists of $M$, embedded in $Z$ with normal bundle  $\kappa^{-1}_M$,
 together with  the image of this surface via the anti-holomorphic map
$\sigma$.
\end{proof}

\begin{example} Let $M={\Bbb CP}_1\times {\Sigma}_{\bf g}$ be the product
of the Riemann sphere with a curve of genus ${\bf g}\geq 2$. Equip
the factors with metrics of curvature $\pm 1$, and let $g$ be the
product metric, which is a scalar-flat K\"ahler metric on $M$.
Because there is only one admissible K\"ahler class on $M$ of
a given volume, $\hat{\cal F} _{\Xi}\equiv 0$ for any
holomorphic vector field $\Xi$. Since $a(M)={\bf sl}(2, {\Bbb C})$,
we therefore have $$h^{2}(Z,{\Theta}_{Z}\otimes {\cal I}_D)=
h^{2}(Z, {\Theta}_{Z,D})=h^{2}(Z,{\Theta}_{Z})=3~ .$$

On the other hand, if $(M, g)$ is  instead the twisted version of the
above example constructed  on the ${\Bbb CP}_1$-bundle associated to any flat
connection on a principal ${\bf SU}(2)$-bundle over ${\Sigma}_{\bf g}$, then,
provided that the given flat connection is
generic in the sense that  its holonomy acts irreducibly on
 ${\bf su} (2)$,  there are no non-trivial
holomorphic vector fields on $M$, and
$$h^{2}(Z,{\Theta}_{Z}\otimes {\cal I}_D)=
h^{2}(Z, {\Theta}_{Z,D})=h^{2}(Z,{\Theta}_{Z})=0~ .$$
\end{example}

\begin{remark}  If  $(M,g)$ is Ricci-flat, the story is utterly different
from that described in Theorem \ref{dt5}. Instead,  one may immediately
read off from Corollary \ref{dt2} that
$$r\equiv 0 \Rightarrow h^2 (Z, {\Theta}_{Z}\otimes {\cal I}_D)=b^+\neq 0.$$
The deformation techniques we are developing here are thus ill-suited to,
say, a K3 surface. Instead, in this hyper-K\"ahler case, when there is
more than one choice of parallel complex structure available,
an unobstructed deformation theory can  be obtained by considering
deformations of $Z$ relative to a fibration over ${\Bbb CP}_1$.
However, in light of the quite definitive theory of Ricci-flat K\"ahler metrics
one obtains from
Yau's solution of the Calabi conjecture \cite{yau2}, there is
little reason  to  pursue this point of view.
\end{remark}

As our first application of this result, let $M$ be a compact K\"{a}hler
surface with a fixed complex structure $J$ and $c^{2}_{1}<0$.  Introduce
${\cal S}$, the moduli space of scalar-flat K\"{a}hler metrics modulo
homothety, and ${\cal A}_M/\br^+$, the projectivized cone of
K\"ahler classes which are  $\cup$-orthogonal to $c_{1}$.  There is a natural
map $\mu:{\cal S}\rightarrow {\cal A}_M/\br^+$ induced by mapping a metric to
its
K\"{a}hler class.

\noindent
\begin{thm}
Let $g$ be a scalar-flat K\"{a}hler metric on
$M$.  Assume that $g$ is not Ricci-flat.

(i) If   $d\hat{\cal F}_{\Xi}|_{[\omega_g ]}\neq 0$
for every non-zero holomorphic vector field $\Xi$ on $M$,
the deformation theories for problems (b), (c), and (d) are all
unobstructed. In particular,
$g$ is a smooth point of the moduli space ${\cal S}$
of problem (b), and ${\cal S}$
has dimension $|\tau (M)|$ near $g$.

(ii) If $M$ carries no non-trivial holomorphic vector fields,
the moduli space ${\cal S}$
is smooth, and   $\mu$ is a
local diffeomorphism between  ${\cal S}$
and ${\cal A}_M/\br^+$.   The set of K\"ahler classes represented by
 scalar-flat  metrics
is therefore  open in the space ${\cal A}_M$ of admissible K\"ahler classes.
 \label{dt6}
\end{thm}

\noindent{\bf Proof}.  The first part is a consequence of the relevant
Kodaira-Spencer
Theorem \ref{ks3}(ii), Corollary \ref{dt2} and Theorem \ref{dt5}.
  The second part follows
from Proposition \ref{dt1} and a simple count of dimensions: from its
definition, ${\cal A}_M/\br^+$ is a manifold of dimension $b_{2}-2$ and this
coincides with the dimension $|\tau(M)|$ by Yau's Theorem \ref{sd3}.
{\hfill \rule{.5em}{1em} \\}

\pagebreak
\section{Ruled Surfaces}
\setcounter{equation}{0}
\subsection{Computing the Futaki Invariant}
\label{foo}

Let $(M, J)$ be a compact complex surface with $c_1^{\Bbb R}(M)\neq 0$,
and suppose that $[\omega ]$ is a K\"ahler
class on $M$ such that the total scalar curvature vanishes---
equivalently, such that $c_1\cup [\omega ] =0$. Then,  by Theorem \ref{sd3},
$M$ must be a ruled surface, which is to say that $(M, J)$ is obtained
from a projectivized rank-2 vector bundle ${\Bbb P}(E)\to \Sigma_{\bf g}$
over a compact complex curve $\Sigma_{\bf g}$ by blowing up $m= |\tau (M)|$
points.

As our eventual goal is to study scalar-flat K\"ahler surfaces,
we will only wish to  consider surfaces $M$ with vanishing
Matsushima-Lichnerowicz
obstruction in the sense of Definition \ref{matlic}. The search is
therefore considerably
narrowed by the following result:

\begin{propn}
Let $M$ be a compact complex surface with an admissible K\"ahler class,
vanishing Matsushima-Lichnerowicz obstruction, and non-trivial
automorphism algebra $a(M)$.
 Suppose also that $M$ is not finitely covered by a
a complex torus.
Then, for some holomorphic line bundle ${\cal L}\to \Sigma_{\bf g}$ over a
compact
complex curve  of genus ${\bf g}\geq 2$,
$M$ is obtained from the minimal
ruled surface
${\Bbb P}({\cal L}\oplus \O )\to \Sigma_{\bf g}$ by blowing up
$|\tau (M)|$ points along the zero section of
 ${\cal L}\subset{\Bbb P}({\cal L}\oplus \O )$. Moreover,
unless $M= \bcp_1\times\Sigma_{\bf g}$,
the space  $a(M)$ of holomorphic vector fields is 1-dimensional,
and is spanned by  the Euler vector
field of ${\cal L}$.\label{lem}
\end{propn}
\begin{proof}
Theorem \ref{sd3} tells us immediately that $M$ is either ruled or covered by
a K3 surface. The latter possibility, however,
 is excluded because
$\Gamma ( {\bf K3}, \Theta )=0$.

 Since $c_1\cdot [\om  ]=0$,
 $c_1^{\Bbb R}$ is  a non-zero primitive
class in $H^{1,1}$, and  $c_1^2 <0$. Thus, with respect to the
 complex orientation, $2\chi +3\tau < 0$.
If the curve $\Sigma_{\bf g}$ has genus $<2$,  we therefore
have an   estimate of the
number of (-1)-curves contained in $M$. Specifically, if
$M$ is obtained by blowing up a  Hirzebruch surface
${\Bbb P}(\O(k)\oplus \O)
\to {\Bbb CP}_1$ at $\ell$ points, we must have $\ell\geq 9$; and if
$M$ is instead obtained by blowing up a minimal ruled surface
$\check{M}\to {\Bbb E}$
over an elliptic curve ${\Bbb E}={\Bbb C}/\Lambda$ at $\ell$ points, then
$\ell>0$.

Let $\Xi\neq 0$ denote a holomorphic vector field on $M$, and let
$\pi : \check{M}\to \Sigma_{\bf g}$ denote a $ {\Bbb CP}_1$-bundle from which
$M$
can be obtained by a blow-up $b: M\to \check{M}$. We then consider the
component of $b_{\ast}\Xi$ normal to the  fibers of $\pi$. Since
the normal bundle of each such fiber is trivial,
 this normal component is constant
up the fibers, so that $(\pi b)_{\ast}\Xi$ is a well-defined
holomorphic vector field on $\Sigma_{\bf g}$. If $\Sigma_{\bf g}$ has
 genus $>1$, this vector field
must vanish.
If, on the other hand, $\Sigma_{\bf g}$ has genus 1, the fact that $b$ involves
blowing up at at least one point forces $b_{\ast}\Xi$, and hence
$(\pi b)_{\ast}\Xi$,
 to have at least
one zero, implying that $(\pi b)_{\ast}\Xi\equiv 0$. Finally, if
$\Sigma_{\bf g}$ has genus 0, and if $(\pi b)_{\ast}\Xi\not\equiv 0$, the
$\ell\geq 9$
blown-up points
of $\pi :F_k\to {\Bbb CP}_1$ must be located on at most 2 fibers of $\pi$;
but then $F_k={\Bbb P}(\O \oplus\O (k))$ admits  sections of $\kappa^{-\ell}$
which vanishes along these two fibers to order $\ell$, and this section lifts
to a non-zero element of $H^0(M, \O (\kappa^{-\ell}))$,
contradicting Corollary \ref{yup}. The vector field  $(\pi b)_{\ast}\Xi$
must therefore
vanish identically, and $ b_{\ast}\Xi$ is tangent to the fibers of
$\pi$. In short $a(M)$ consists strictly of {\em vertical} vector fields.

 By the Matsushima-Lichnerowicz assumption, the identity
component of the automorphism group of $M$ is the complexification of a
compact group. We therefore have a non-trivial holomorphic  vector field
$\Xi$    on $M$ whose
imaginary part  $\xi$ generates an $S^1$-action, and which itself generates
a $\bc^{\times}$-action; for brevity's sake, we shall henceforth
refer to any such  $\Xi$  as a {\em  periodic holomorphic vector field}.
 On the other hand,
the minimal model
 $\pi : \check{M}\to \Sigma_{\bf g}$ may be represented in the form
${\Bbb P}(E)\to \Sigma_{\bf g}$ for a rank 2 holomorphic vector bundle
$E\to \Sigma_{\bf g}$
which is completely
specified once an arbitrary  line bundle  $\wedge^2E$
is chosen, subject to the condition $c_1(E)\equiv w_2(\pi)\bmod 2$.
The vector field  $ b_{\ast}\Xi$ is then uniquely specified by a
trace-free holomorphic section $A$ of ${\cal E}nd (E)$. The determinant
of $A$ is  a holomorphic function on $\Sigma_{\bf g}$, hence a constant.
On the other hand, since  $\Xi$ is periodic, $A$
is diagonalizable, and  $A$ must be a half-integer  multiple of
 $$\left[\begin{array}{cc}1&0\\0&-1
\end{array} \right]~ .$$
The vector bundle $E$ thus globally splits as a direct sum of the eigenspaces
of
 $A$, and, twisting by a line bundle,
we may therefore take $E={\cal L}\oplus \O$, so that  $\Xi$ becomes a
constant  multiple
of the Euler vector field on ${\cal L}$. We henceforth normalize this constant
to be 1.
The blown-up points must all occur at zeroes of $\Xi$, namely either
at the zero section of ${\cal L}$ or at the ``infinity section''
corresponding to the $\O$ factor. The latter possibility may be
reduced to the former by noticing that the proper transform of a fiber
 through exactly one blown-up
point is a (-1)-curve, which may therefore be blown down, thereby leading to
a different
minimal model. In our case, iteration of this procedure
allows us to replace  blown-up
points ``at the infinity section''  by  blown-up points
``at the zero section,'' at the small price of
twisting our line bundle $\cal L$ by the divisor of the relevant
points of $\Sigma_{\bf g}$.

The space of vertical vector fields on the minimal model $\check{M}$ is
now precisely $\Gamma (\Sigma_{\bf g}, \O \oplus {\cal L}\oplus {\cal L}\* )$,
with the Lie algebra structure induced by identifying
$(u,v,w)\in\Gamma (\Sigma_{\bf g}, \O \oplus {\cal L}\oplus {\cal L}\* )$
with the matrix
$$\left[\begin{array}{cc} u&v\\w&-u\end{array}\right]~ .$$
If $M=\check{M}$, this algebra is itself required to be reductive,
implying $H^0(\Sigma_{\bf g}, {\cal L})\neq 0 \Leftrightarrow
H^0(\Sigma_{\bf g}, {\cal L}\* )\neq 0$; we conclude that
either   ${\cal L}$ is trivial and $M=\bcp_1\times \Sigma_{\bf g}$, or
else
$a(M)$ is 1-dimensional.
If, on the other hand,
$M$ is obtained by
 blowing up points on the zero section of $\cal L$, the vector field
$(u, v, w)$
 lifts to $M$ iff $v=0$, and
$a(M)=\{ (u,0,w )\}$; thus, if $\tau (M)\neq 0$,
$a(M)$ is reductive iff $\Gamma (\Sigma_{\bf g},{\cal L}\* )=0$.
Thus, provided that $M\neq
\bcp_1\times \Sigma_{\bf g}$,
$a(M)$ is 1-dimensional, with the Euler vector
fields $\Xi$ (corresponding to $(u,v,w)=(\frac{1}{2},0,0)$) as a basis.

Finally, we observe that $M$ must have genus $\geq 2$. Indeed, the Euler
field $\Xi$ is a vector field on $\check{M}$ which vanishes at all the
points which are to be blown up. Let $\ell$ be the greatest multiplicity
with which any point is to be blown up, and, assuming ${\bf g}=0,1$,
 let $\Upsilon\not\equiv 0$ be any vector field on $\Sigma_{\bf g}$.
Then $(\Xi\wedge\Upsilon)^{\otimes \ell }$ lifts to $M$ as a non-trivial
section of $\kappa^{-\ell }$, contradicting Proposition \ref{yup}.
Hence ${\bf g}\geq 2$.
\end{proof}

Our  goal is now to calculate the
Futaki invariant of $(M, [\omega ], \Xi  )$, where $M$ is
in normal form described in the above Proposition,
$[\om ]$ is an admissible K\"ahler class and  $\Xi$ is the
Euler vector field.  We proceed by a symplectic
quotient construction in the spirit of \cite{L}.

The invariant we seek to compute
is known \cite{cal2}
to be independent of the representative $\omega\in [\omega ]$,
so we may assume (by averaging) that  $\omega$ is invariant under
the $S^1$-action generated by $\xi=\Im \Xi$.
Since $0={\pounds}_{\xi } \omega = d (\xi\rfloor\omega )$,
we see that  $ \nu : =\omega (  \xi , \cdot )$
is closed; and,  on the other hand,
 any real harmonic 1-form on our compact K\"ahler manifold
is the real part  of a holomorphic 1-form, and so
must everywhere be orthogonal to $\xi =\Im \Xi $, as may either
be seen directly from the our explicit description of $(M,\Xi  )$, or
deduced as a consequence
of the  maximum
principle for pluriharmonic functions and the fact that $\Xi $ has zeroes.
Harmonic theory therefore yields
$\nu := (dd^{\ast}+d^{\ast}d)G\nu=dd^{\ast}G\nu$,
where $G$ is the Green's operator, and $\Re\Xi  := \mbox{grad} f$
for a unique function $f:= d^{\ast}G\nu$,
called the   holomorphy potential \cite{besse}\cite{lich} of $\Xi $,
such that $\int_M f\vol=0$. As we saw in \S \ref{key},
the Futaki invariant ${\cal F}(\Xi ,  [\omega ])$ of $(M,[ \omega ])$
is then given  by
$${\cal F}(\Xi ,  [\omega ])= -{\textstyle \frac{1}{2}} \int_M fs\vol~ . $$

The symplectic vector field
$\xi = \Im \Xi $ generating the $S^1$-action is now a globally Hamiltonian
vector field, meaning that
$$\omega (  \xi , \cdot ) = dt$$
for a smooth (``Hamiltonian'')
function $t:M\to {\Bbb R}$; indeed, anything of the form
of the form $t=f+c$ will do. We could, of course,
 choose our constant
$c$ to vanish, but we will instead find it convenient to
choose $c$  so
that $\max t=-\min t=a$, and $t: M \mapsonto [-a,a]$.
(As we shall see in a moment, the intrinsic significance of the
number $a$ is that  $\int_F [\omega ]=4\pi a$, where $F$ is
any fiber of $M\to \Sigma_{\bf g}$.)
Fortunately, because we have assumed that $c_1\cup [\omega ]=0$,
this will not interfere with our calculation of the Futaki invariant
because
\bea\int_Mts\vol&=&\int_M (f+c)s\vol\\&=&
\int_M fs\vol+c\int_M s\vol\\&=&
\int_M fs\vol\\&=&
-2{\cal F}(\Xi ,  [\omega ])~ .\eea

The only isolated critical points of $t$ occur at those zeroes of $\Xi$
which occur at the intersection of an exceptional curve and the
proper transform of a fiber; since such a fixed point is attractive along the
exceptional curve and repulsive along the proper transform of a fiber,
such a critical point has index 2. On the other hand,
 the maxima and minima of $t$ occur along a pair of
holomorphic curves,  $C_0=t^{-1}(-a)$ and $C_{\infty}=t^{-1}(a)$,
which are just  the proper transforms of the
 ``zero'' and ``infinity'' sections of
$ {\Bbb P}({\cal O}\oplus {\cal L})\to \Sigma_{\bf g}$. We now have
 a projection
$\pi : M\to \Sigma_{\bf g}\times [-a,a]$ given by the product of
the ruling $M\to \Sigma_{\bf g}$
and the Hamiltonian $t$. Because the ${\Bbb C}^{\times}$-action preserves
the ruling, every fiber of $\pi$ consists of exactly one orbit of the
$S^1$-action. (This is really a consequence \cite{atiyah}\cite{mum}
of the fact that
$\Sigma_{\bf g}$ is both the symplectic and stable quotient
of $M$ by the ${\Bbb C}^{\times}$-action.) Let
$q_1 , \ldots , q_m\in \Sigma_{\bf g}\times (-a,a)$
be the images of the isolated fixed points of the action, and let
$X :=   [\Sigma_{\bf g}\times (-a,a)]-\{ q_1 , \ldots , q_m\}$
denote the  set of regular values of $\pi$.  If $Y\subset M$ is the set of
regular points, then $\pi : Y\to X$ is a principal $S^1$-bundle, and,
by taking the orthogonal complement of the $S^1$ orbits with respect to
the K\"ahler metric, we endow   $Y\to X$ with a connection form
$\theta$. If $z=x+iy$ is any complex local coordinate on $\Sigma_{\bf g}$,
we may then express the given  K\"ahler metric $g$ on $Y\subset M$
in the form\be
g=vw(dx^{\otimes 2}+dy^{\otimes 2})+w~dt^{\otimes 2}+w^{-1}\theta^{\otimes 2}~,
\label{met}\ee
for positive functions $v, w>0$ on $X$, while the
complex structure $J$ is given by
\bea
dx&\mapsto &dy\\dt&\mapsto& w^{-1}\theta
\eea
so that the K\"ahler form is given by
$$\omega = dt\wedge \theta +vw~ dx\wedge dy~ .$$
Since the complex structure $J$ is integrable, the differential ideal
$${\cal J} = \langle dx + i dy ,  w dt + i\theta\rangle$$ must satisfy
$d{\cal J}\subset {\cal J}$; explicitly, this means that
\bea d (w dt + i\theta)&=& dw\wedge dt + id\theta \\
&=&\varphi\wedge (dx +i dy)\eea
for some complex-valued 1-form $\varphi$ on
$X$, and, because $\theta$ and $w$ are real,
 this is in turn equivalent to
\be d\theta \equiv w_x dy\wedge dt+
w_y dt\wedge dx \bmod dx\wedge dy~ .\label{ka2} \ee
The K\"ahler condition $d\omega=0$  now reads
\be 0=d(dt\wedge \theta +vw dx\wedge dy)= -dt\wedge d\theta +
(vw)_t dt\wedge dx\wedge dy~,
\label{ka1}\ee
so that  the curvature of our
 $S^1$-connection $\theta$ is now completely determined by $v$ and $w$:
$$d\theta = w_x dy\wedge dt+ w_y dt\wedge dx + (vw)_t dx\wedge dy~ .$$
In particular, we conclude that
$w_{xx}+w_{yy}+(vw)_{tt}=0$.

Notice that  equation (\ref{met}) says that the metric on any fiber
$F$ of $M\to \Sigma_{\bf g}$ is given by
$$g|_F=wdt^2+ w^{-1}d\vartheta^2~,$$
where  $\vartheta\in [0,2\pi ]$ is a fiber coordinate in a gauge
chosen such that the connection form $\theta$ has no $dt$ component; the
area form on $F$ is just\footnote{This generalizes to a  simple
relationship between volumes and moment maps for torus actions that
is sometimes called the
``Archimedes Principle'' \cite{archie}\cite{atiyah}.}
$$\omega |_F= dt\wedge d\vartheta ~,$$
and the area of $F$ is therefore $4\pi a$. On the other hand, since this metric
is  smooth at the ``south pole'' $t=-a$ of the 2-sphere $F$,
letting $r$ denote the Riemannian distance from
the south pole, we have
$$w~dt^2+ w^{-1}~d\vartheta^2=dr^2+(r^2+O(r^4))~d\vartheta^2~,$$
so that  $dt=r(1+O(r^2))~dr$, $t+a=\frac{r^2}{2}+O(r^4)$, and
$w^{-1}=2(t+a)+O((t+a)^2)$. Similarly,
$w^{-1}=-2(t-a)+O((t-a)^2)$ near $t=a$.
Thus $\ell =
w^{-1}$, which  is a smooth function on $M$ because it
 represents the
square of the length of the Killing field $\xi$,
descends to a differentiable
function on $\Sigma_{\bf g}\times [-a,a]$ which
vanishes at the boundary and
satisfies $\frac{d\ell }{dt}=\mp 2$ at $t=\pm a$.
 At the same time,
equation (\ref{met}) tells us that $\lim_{t\to -a} vw~ dx\wedge dy=
\omega |_{C_0}$, while   $\lim_{t\to a} vw~ dx\wedge dy=
\omega |_{C_{\infty}}$. Thus $v$ is smooth up to the boundary of
$\Sigma_{\bf g}\times [-a,a]$, and moreover
\bea \left. v\right|_{t=\pm a}&=&0 \\
\left. v_t dx\wedge dy\right|_{t=-a}  &=&
  \left. \hphantom{-}2\omega \right|_{C_0}  \\
\left. v_t dx\wedge dy\right|_{t=a\hphantom{-}}  &=&
   \left. -2\omega \right|_{C_{\infty }}~ .  \eea

Since $\Xi $ is a holomorphic vector field and the $(2,0)$-form
$\mu := dz\wedge (w~dt + i\theta )$ has the property that
$\Xi  \rfloor \mu = 2 ~dz$ is a holomorphic form, $\mu$ must itself
be holomorphic; thus $\frac{1}{4}
\mu\wedge\overline{\mu }= w dx\wedge dy\wedge dt \wedge
\theta$ is the  volume form of a holomorphic frame.
On the other hand, the metric volume form is
$${\textstyle \frac{1}{2}}\omega\wedge\omega=vw~ dx\wedge dy\wedge dt\wedge
\theta~,
$$
so that the Ricci
form of $g$ must be
$$\rho = -i\partial \overline{\partial }\log \left(
\frac{vw~ dx\wedge dy\wedge dt\wedge \theta}{w~ dx\wedge dy\wedge dt\wedge
\theta}
\right)= -i\partial \overline{\partial }\log v~ .$$
The scalar curvature density of $g$ thus is given  in terms of $u:=\log v$ by
\bea s\vol&=& {\textstyle \frac{1}{2}} s~\omega\wedge\omega\\
&=&2  \omega\wedge\rho\\
&=& -2i\omega\wedge \partial \overline{\partial }u \\
&=& -2i\omega\wedge [-{\textstyle\frac{i}{2}} dJd u ]\\
&=& -\omega\wedge dJd u \\
&=& -\omega\wedge d[u_xdy- u_ydx+u_tw^{-1}\theta ]\\
&=& -[dt\wedge\theta + vw~dx\wedge dy ]\wedge
d[u_xdy- u_ydx+u_tw^{-1}\theta ]\\
&=& [u_{xx}+ u_{yy}+ vw(w^{-1}u_t)_t]~dx\wedge dy\wedge
dt \wedge\theta + (u )_tw^{-1}d\theta \wedge dt \wedge \theta \\
&=& [u_{xx}+ u_{yy}+ vw(w^{-1}u_t)_t +w^{-1}u_t (vw)_t
]~dx\wedge dy\wedge
dt \wedge\theta  \\
&=& [u_{xx}+ u_{yy}+ (vu_t)_t
]~dx\wedge dy\wedge
dt \wedge\theta  \\
&=&  [( \log v)_{xx}+ ( \log v)_{yy}+ v_{tt}]~dx\wedge dy\wedge dt \wedge\theta
{}~ .
\eea
In particular,
\be s=  \frac{( \log v)_{xx}+ ( \log v)_{yy}+ v_{tt}}{vw} ~ . \label{scal}
\ee

Let us now rephrase the above results in more global terms. Because
$g$, $w$ and $dt$ are globally defined, it follows from  equation (\ref{met})
that, for $-a< t < a$,
 $$g^v(t):=v~(dx^2+dy^2)$$
 is a well-defined $t$-dependent K\"ahler metric
on $\Sigma_{\bf g}$, with K\"ahler form $$\omega^v (t):= v~dx\wedge dy~ .$$
Moreover,
\bea \left.\omega^v\right|_{t=\pm a}&=&0 \\
\left.\frac{d }{d t} \omega^v \right|_{t=-a} &=&
   \left.\hphantom{-}2\omega \right|_{C_0}  \\
\left.\frac{d }{d t} \omega^v \right|_{t=a\hphantom{-}} &=&
 \left. - 2\omega \right|_{C_{\infty }}~ .  \eea
The Ricci form of this metric is
\bea \rho^v(t)&=& -i(\frac{\partial^2}{\partial z\partial
 \overline{z}} \log v )dz\wedge d\overline{z}\\ &=&
-{\textstyle \frac{1}{2}}[(\log v )_{xx}+ (\log v )_{yy}]~dx\wedge dy ~,
\eea
 so that our formula for the density of scalar curvature may be written
globally on $Y\subset M$ as
$$ s\vol = [-2\rho^v + \frac{d^2}{dt^2} \omega^v]\wedge
dt\wedge \theta~ .$$
Hence
\bea \int_M ts\vol &=& \int_Y ts\vol
\\&=&
\int_Y t[-2\rho^v + \frac{d^2}{dt^2} \omega^v]\wedge dt\wedge \theta
\\&=&
2\pi \int_{-a}^at[\int_{\Sigma_{\bf g}}
-2\rho^v + \frac{d^2}{dt^2} \omega^v]~dt
\\&=&
2\pi \int_{-a}^at\frac{d^2}{dt^2}[\int_{\Sigma_{\bf g}}\omega^v]~dt
{}~~+2\pi \int_{-a}^at [\int_{\Sigma_{\bf g}}
-2\rho^v  ]~dt
\\&=&
2\pi \int_{-a}^at\frac{d^2}{dt^2}[\int_{\Sigma_{\bf g}}\omega^v]~dt
-8\pi^2 \chi (\Sigma_{\bf g})\int_{-a}^at  ~dt
\\&=&
2\pi \left[t\frac{d}{dt}\int_{\Sigma_{\bf g}}\omega^v\right]^a_{-a}-
2\pi \int_{-a}^a\frac{d}{dt}[\int_{\Sigma_{\bf g}}\omega^v]~dt
\\&=&
2\pi \left[t\frac{d}{dt}\int_{\Sigma_{\bf g}}\omega^v\right]^a_{-a}-
2\pi \left[\int_{\Sigma_{\bf g}}\omega^v\right]^a_{-a}
\\&=&
2\pi \left[t\frac{d}{dt}\int_{\Sigma_{\bf g}}\omega^v\right]^a_{-a}
\\&=&
2\pi a\left[\int_{\Sigma_{\bf g}} \left.\frac{d}{dt}\omega^v\right|_{t=a}
 +
\int_{\Sigma_{\bf g}} \left.\frac{d}{dt}\omega^v\right|_{t=-a}\right]
\\&=&
2\pi a\left[\int_{C_0}2\omega -
\int_{C_{\infty}}2\omega \right] =
4\pi a\left[\int_{C_0}\omega -
\int_{C_{\infty}}\omega \right]
\\&=&
\left[\int_{C_0}\omega -
\int_{C_{\infty}}\omega \right]\int_F\omega
\eea
In conclusion, we have
\bea {\cal F}(\Xi ,  [\omega ])&=&
-{\textstyle \frac{1}{2}}\int_M fs\vol
\\&=& -{\textstyle \frac{1}{2}}\int_M ts\vol\\&=&
-{\textstyle \frac{1}{2}}\left[\int_{C_0}\omega -
\int_{C_{\infty}}
\omega\right]\int_F\omega \\&=&
{\textstyle \frac{1}{2}}\left[\int_{C_{\infty}}
\omega -\int_{C_0}\omega\right]\int_F\omega
\eea
where $F$ is any fiber of $M\to \Sigma_{\bf g}$. We have thus
proved the following:

\begin{thm} Let $M$ be any compact complex surface equipped with
 an admissible K\"ahler class $[\omega ]$
and
a holomorphic ${\Bbb C}^{\times}$-action which is free on an open
dense set. Assume that $c_1^{\Bbb R}(M)\neq 0$, and let
$\Xi \in \Gamma (M, \O
(TM))$  denote the holomorphic vector field which generates the action.  Thus
$M$ is a ruled surface $M\to \Sigma_{\bf g}$ of genus  ${\bf g}\geq 2$,
the generic fiber $F$ of which is the closure of an orbit,
while the ``attractive'' and ``repulsive'' fixed curves
 $C_0$ and $C_{\infty}$ of the action are sections of the
projection $M\to \Sigma_{\bf g}$.
The Futaki invariant of $(M, J, [\omega ], \Xi  )$ is then given by
$${\cal F}(\Xi ,  [\omega ])=
{\textstyle \frac{1}{2}}\left[\int_{C_{\infty}}
\omega -\int_{C_0}\omega\right]\int_F\omega ~ .$$ \label{local}
\end{thm} \hfill \rule{.5em}{1em}

\begin{remark} We have computed the Futaki invariant only for a single
 vector field $\Xi$. However, if the Matsushima-Lichnerowicz obstruction
vanishes, Proposition \ref{lem} tells us that either $a(M)$ is spanned by
$\Xi$ or else $M=\bcp_1\times \Sigma_{\bf g}$.
In the former case, the ${\Bbb C}$-linearity of ${\cal F}(\cdot, [\omega ])$
tells us that our computation completely determines the
${\cal F}(\cdot, [\omega ])$.
In the exceptional  case $M=\bcp_1\times \Sigma_{\bf g}$, the
Futaki character
vanishes, since the product of two constant curvature metrics has constant
scalar curvature.
\end{remark}

Use of the fact that $c_1\cdot [\omega ]=0$ allows one to rewrite
the Futaki invariant in  interesting equivalent ways. As previously
indicated, we will always normalize
the minimal model of our ruled surface with holomorphic vector field
by putting it in the form
${\Bbb P}({\cal L}\oplus {\cal O})\to \Sigma_{\bf g}$, where
all the blown-up points {\em are  on the zero section
of} ${\cal L}\hookrightarrow {\Bbb P}({\cal L}\oplus {\cal O}):
\zeta\mapsto [\zeta ,1]$, and so correspond to fibers of
${\cal O}\subset {\cal L}\oplus {\cal O}$. For simplicity, let us assume
for the moment that the blown-up points in ${\Bbb P}({\cal L}\oplus {\cal O})$
are all distinct, and so give rise to $m=|\tau (M)|$
 exceptional rational curves $E_j$ of self-intersection $-1$.
Let us associate $m$ ``weights'' $w_j>0$, $j=1, \ldots , m$,
to the K\"ahler class $[\omega ]$ by defining
$$w_j:=\frac{\int_{E_j}\omega }{\int_{F}\omega }~ .$$
The homology classes of
$C_{\infty}, F, E_1, \ldots , E_m$ form a basis for $H_2 (M)$, and the
intersection form of $M$ is
$$
\left[\begin {array}{ccccc}
-k&1&&&\\1&0&&\\&&-1&&\\&&&\ddots&\\&&&&-1
\end{array}\right]$$
with respect to this basis, where $k:=\mbox{deg}({\cal L})$; it follows that
the Poincar\'e dual of $[\omega ]$ can be expressed in this basis as
$A (1, B+k, -w_1, \ldots , -w_m)$, where $A=\int_F\omega $ and
$AB=\int_{C_{\infty }}\omega$. Since $C_{\infty }$ has genus $\bf g$
and self-intersection $-k$, whereas $F$, $E_1, \ldots , E_m$ have
genus 0 and self-intersection $0, -1, \ldots , -1$, the adjunction
formula may be used to rewrite the
condition
$c_1\cdot [\omega ]=0$ as
$$0=[2(1-{\bf g})-k, 2, 1, \ldots , 1]
\left[ \begin{array}{c}1\\B+k\\ -w_1\\ \vdots \\-w_m\end{array}\right]
=2(1-{\bf g})+k+2B-\sum_{j=1}^mw_j~,$$
so that $B=[-k+2({\bf g}-1)+\sum_{j=1}^mw_j ]/2$, and
$$\int_{C_{\infty }}\omega = \frac{A}{2}[-k+2({\bf g}-1)+\sum_{j=1}^mw_j ]~ .$$
On the other hand, the picture is symmetrical between ${C_{\infty }}$
and $C_0$ as long
as we remember to replace $k=-C_{\infty }\cdot C_{\infty }$ with
$m-k =-C_{0}\cdot C_{0}$ and replace the $E_j$ with new exceptional curves
$\hat{E}_j$ such that $[E_j]+[\hat{E}_j]=[F]$, resulting in a
replacement of the  weights $w_j$ by new weights  $1-w_j$; the upshot is
that
\bea\int_{C_{0 }}\omega &=& \frac{A}{2}[k-m+2({\bf g}-1)+\sum_{j=1}^m(1-w_j) ]
\\&=&\frac{A}{2}[k+2({\bf g}-1)-\sum_{j=1}^mw_j ]~ .\eea
The Futaki invariant can therefore be rewritten as
\bea {\cal F}(\Xi ,  [\omega ])&=&
{\textstyle \frac{1}{2}}\left[\int_{C_{\infty}}
\omega -\int_{C_0}\omega\right]\int_F\omega \\&=&
-\frac{A^2}{4}\left[ (k+2({\bf g}-1)-\sum_{j=1}^mw_j)-
(-k+2({\bf g}-1)+\sum_{j=1}^mw_j)\right]\\&=&
-\frac{A^2}{2} [k-\sum_{j=1}^mw_j]\\&=&
-\frac{A^2}{2} [\mbox{deg}({\cal L})-\sum_{j=1}^mw_j]~ .
\eea
We therefore have the following :
\begin{propn} Let $M$ be obtained from the minimal model \linebreak
${\Bbb P}({\cal L}\oplus {\cal O})\to \Sigma_{\bf g}$, ${\bf g}\geq 2$,
by blowing up $m$  points on the zero section of $\cal L$.
Let $[\omega ]$ be a K\"ahler class on $M$ such that
 $c_1\cup [\omega ]=0$, and let the weights $w_j\in ~]0,1[$ be defined by
$\int_{E_j}\omega =w_j\int_F\omega$, where $F$ is a typical fiber of
$M\to  \Sigma_{\bf g}$ and the $E_j$, $j=1,\ldots ,m$
 are the exceptional
curves corresponding to the blown-up points. Let $\Xi $ denote the
vector field on $M$ corresponding to the Euler field on $\cal L$. Then
$${\cal F}(\Xi ,  [\omega ])=0~~\Longleftrightarrow
{}~~\sum_{j=1}^mw_j=\mbox{deg}({\cal L})~ .$$\label{foot}
\end{propn}
\begin{proof}
We assumed for simplicity in the previous discussion that $M$ is obtained from
its minimal model by blowing up {\em distinct} points. However, the argument
  goes through without change
provided that, in defining the weights $w_j$, one replaces the
integrals of $\omega$ over the exceptional divisors $E_j$ with
integrals of $\omega$ over the corresponding  homology classes,
each of which can be represented by a chain of rational curves with
intersection matrix
$$\left[ \begin{array}{cccc}
-2\hphantom{-}&1&&\\1&\ddots&1&\\
&1&-2\hphantom{-}&1\\&&1&-1\hphantom{-}\end{array}\right]~ .$$
\end{proof}
\begin{cor} Let $(M, \Xi)$ be as in Proposition \ref{foot}.
Then $M$ carries an admissible K\"ahler class $[\omega ]$
for which
${\cal F}( \Xi , [\om ])=0$  iff one of the  following
holds:
\begin{description}
\item{(a)}  $0=\mbox{deg}({\cal L})=m$; or
\item{(b)}   $0<\mbox{deg}({\cal L})<m$.
\end{description}\label{tri}
 \end{cor}
\begin{proof} The necessity of these conditions is an immediate consequence of
Proposition \ref{foot}. For sufficiency, one may either invoke
 Proposition \ref{nakai},  or else  wait for the explicit construction
in \S \ref{next} below.
\end{proof}

\begin{remark} Notice that $m=1$ is
excluded, since $(b)\Rightarrow m\geq 2$. Also notice that
the vanishing of the Futaki invariant implies the vanishing of
the Matsushima-Lichnerowicz obstruction.
Indeed, if $m=0$ we either have $\cal L$ is trivial, and
$M=\bcp_1\times \Sigma_{\bf g}$, or else
$\Gamma (\Sigma_{\bf g}, {\cal L})=
\Gamma (\Sigma_{\bf g}, {\cal L}\* )=0$, so that $a(M)$
is generated by the Euler field $\Xi$;  if $m>0$,
 $\Gamma (\Sigma_{\bf g}, {\cal L}\* )=0$ and again
$a(M)$
is generated by the Euler field $\Xi$.
\end{remark}
\begin{cor} Suppose that $(M, [\omega ])$  is a compact complex surface with
admissible K\"ahler class and vanishing Matsushima-Lichnerowicz obstruction.
Let $\Xi\not\equiv 0$ be any non-trivial holomorphic vector field, and
consider the
 restricted Futaki functional $\hat{\cal F}_{\Xi}$ of
Definition \ref{two}. Then
 $$d\hat{\cal F}_{\Xi}|_{[\omega ]}=0 ~\Longleftrightarrow ~ \tau (M)=0.$$
\label{van}
\end{cor}
\begin{proof} By Proposition \ref{lem}, the only non-ruled surfaces we need
consider are tori and their (hyper-elliptic) quotients, for which
both $\tau$ and $\cal F$ vanish. For the ruled surfaces, the result
follows immediately from Theorem \ref{local}.
\end{proof}

Finally, as an aside, let us  observe
that  the condition ${\cal F}(\Xi ,  [\omega ])=0$
can now be restated in terms  of parabolic stability
in the sense of Seshadri \cite{sesh}\cite{mehta}.
We consider the vector bundle ${\cal V}={\cal L}\oplus \O$ of our minimal
model,
equipped with  1-dimensional subspaces
 $L_j$ in some fibers of ${\cal V}$ which represent the exceptional divisors
$E_1, \ldots , E_m\subset M$. (Thus the
subspaces $L_j$ are contained the  $\O$ factor of ${\cal L}\oplus \O$.)
Let $w_1, \ldots , w_m$ denote the weights as before, and let
$\alpha_j < \beta_j$ be arbitrary numbers in
 $[0,1]$ such that $w_j=\beta_j-
\alpha_j$. The criterion
$$\sum_{j=1}^mw_j=\mbox{deg}({\cal L})~,$$
is precisely equivalent to the statement that
$({\cal V}, \{( L_j,\alpha_j, \beta_j)\})$
is parabolically
quasi-stable\footnote{
With this choice of terminology, stable $\Rightarrow$ quasi-stable
$\Rightarrow$
semi-stable.}, in
the sense that, for every line sub-bundle $L\subset {\cal V}$, we have
$$\mbox{pardeg}(L)\leq \frac{1}{2}\mbox{pardeg}({\cal V}),$$
with equality iff $L$ is a {\em direct summand} of ${\cal V}$;
here the {\em parabolic degree} of a line sub-bundle $L\subset {\cal V}$
is defined to be
$$\mbox{pardeg}(L):=\mbox{deg}(L)+\sum_{\{j|L_j\not\subset L\}}\alpha_j
+\sum_{\{j|L_j\subset L\}}\beta_j~,$$
whereas
$$\mbox{pardeg}({\cal V}):=\mbox{deg}({\cal V})+\sum_{j=0}^m\alpha_j
+\sum_{j=0}^m\beta_j~ .$$

\begin{cor} The vanishing of the
Futaki invariant for a non-trivial
${\Bbb C}^{\times }$-action  on a
ruled surface with  admissible  K\"ahler class is
equivalent to the quasi-stability of the
the corresponding parabolic bundle.
\end{cor}

\subsection{Scalar-Flat Ruled Surfaces with Vector Fields}
\label{next}
In the last section, we analyzed the Futaki invariant
of compact complex surfaces
with  periodic holomorphic vector fields which admit admissible K\"ahler
classes.
(Recall that {\em admissible} means that the total
scalar curvature of a metric in the class vanishes.) Since the
Matsushima-Lichnerowicz Theorem and the Futaki invariant
are obstructions to the existence of  constant scalar curvature
 K\"ahler metrics
in the given class, this  gives us a rough classification
of those surfaces
with holomorphic vector fields which might admit  scalar flat K\"ahler metrics.
In this section we will review a construction \cite{L2} \cite{L5} of compact
scalar-flat K\"ahler
surfaces, and use it to observe that this ``rough''
 classification is in fact perfectly sharp.

The idea is to reverse the symplectic quotient construction of the
last section. Let $\Sigma_{\bf g}$ be any compact complex curve
of genus  $\geq 2$, and let $h_{\Sigma}$ be the unique Hermitian metric
on $\Sigma_{\bf g}$ of constant curvature $-1$. We can then give the 3-manifold
$\Sigma_{\bf g}\times (-1,1)$ a hyperbolic structure
by introducing the constant curvature $-1$ metric
$$h:=\frac{h_{\Sigma}}{(1-t^2)}+\frac{dt^2}{(1-t^2)^2}~ .$$
Let $q_1, \ldots , q_m\in \Sigma_{\bf g}\times (-1,1)$ be arbitrary points,
and associate to each
the Green's function $G_j$, defined by
$$\Delta G_j=2\pi \delta_{q_j}, ~~~\lim_{t\to \pm 1} G_j = 0 ~,$$
where $\Delta =-\star d\star d$ is the
(positive) Laplace-Beltrami operator of $h$.
We define $V:=1+\sum_{j=1}^mG_j$, so that
$$\Delta V=2\pi \sum_{j=1}^m\delta_{q_j}, ~~~\lim_{t\to \pm 1} V=1 ~ .$$
 $V$ extends smoothly to  $(\Sigma_{\bf g}\times [-1,1])-\{
q_1, \ldots , q_m\}$ and satisfies $V\geq 1$.

On $[\Sigma_{\bf g}\times (-1,1)]-\{
q_1, \ldots , q_m\}$, the 2-form
$$\alpha := \frac{1}{2\pi }\star dV$$
is now
 closed, and its integral on a small sphere around any one of the
$q_j$'s is $-1$. On the other hand, if $1-\epsilon > \max_{j} t(q_j)$, then
one may check \cite{L2} that
\be\int_{\Sigma_{\bf g}\times \{ 1-\epsilon\}}\alpha = -\sum_{j=1}^m\left(
\frac{1+t(q_j )}{2}\right) ~,\label{int}\ee
so that, setting $w_j:=[1+t( q_j)]/2$, we have
\be [\frac{1}{2\pi }\star dV]\in H^2_d
\left( [\Sigma_{\bf g}\times (-1,1)]-\{ q_1, \ldots , q_m\},
{\Bbb Z} \right) ~~\Longleftrightarrow ~~\sum_{j=1}^mw_j\in {\Bbb Z}~,
\label{foutre}\ee
since the second homology of $[\Sigma_{\bf g}\times (-1,1)]-\{
q_1, \ldots , q_m\}$ is generated by the homology
classes of $\Sigma_{\bf g}\times \{ 1-\epsilon\}$ and
 $m$ small spheres centered at the punctures
$q_1, \ldots , q_m$.
If we assume this condition is met,  the Chern-Weil theorem guarantees that
we can then
find a principal $S^1$-bundle  $$\pi_0: M_0\to [\Sigma_{\bf g}\times (-1,1)]-\{
q_1, \ldots , q_m\}$$ with a connection 1-form $\theta$ for which the
curvature is
$$d\theta = \star dV~ .$$
Notice that, even modulo gauge equivalence, the pair $(M_0, \theta )$ is
by no means unique, since our base is not simply connected; instead,
the group $H^1(  \Sigma_{\bf g}, S^1)$ of flat circle bundles on
$\Sigma_{\bf g}$ acts freely and transitively on the orbits by tensor
product. Given a choice of $(M_0, \theta )$ we then
equip $M_0$ with the Riemannian metric
$$g:= (1-t^2)[Vh+V^{-1}\theta^2]~ . $$
If we identify the universal cover of $\Sigma_{\bf g}$  with the upper
half-plane $y=\Im z >0$ in $\Bbb C$, the metric can be more explicitly
 written in the
form
\bea g &=& (1-t^2)\left[ V \frac{dx^2+dy^2}{y^2(1-t^2)}+V\frac{dt^2}{(1-t^2)^2}
+V^{-1}\theta^2\right]
\\ &=&vw~(dx^2+dy^2)+w~dt^2+w^{-1}\theta^2, \eea
where
$$w=\frac{V}{1-t^2}$$
and
$$v=\frac{1-t^2}{y^2}~ .$$
Since the equation $d\theta =\star dV$ can now be rewritten
as
$$d\theta =w_x~dy\wedge dt+w_y~dt\wedge dx+(vw)_t~dx\wedge dy~,$$
our  calculations (\ref{ka2}) and (\ref{ka1}) show that $g$
is K\"ahler with respect to the integrable complex structure
\bea
dx&\mapsto &dy\\dt&\mapsto& w^{-1}\theta
\eea
Moreover, since
$$(\log v)_{xx}+(\log v)_{yy}+v_{tt}=0~,$$
we conclude from equation (\ref{scal}) that $g$ is scalar-flat.

We can now compactify $M_0$ by adding two copies of $\Sigma_{\bf g}$,
corresponding to $t=\pm 1$,
and $m$ isolated points, corresponding to $q_1, \ldots , q_m$. This
 compactification $M$ can then \cite{L2} be made into a smooth manifold
in such a way that the metric $g$ and the complex structure $J$ extend to
$M$, giving us a compact scalar-flat K\"ahler surface $(M,g)$.
The bundle projection $\pi_0$ now extends to a smooth map
$$\pi: M \to \Sigma_{\bf g}\times [-1,1]~,$$
and the original $S^1$-action extends to  an action on $M$ for which
$\pi$ is projection to the orbit space; the points  added to $M_0$ in order
to obtain $M$ are
precisely the fixed points of the action.
The tautological   projection
 $\mbox{pr}_1\pi: M\to  \Sigma_{\bf g}$ induced by $\pi$ and
the first-factor projection
$\mbox{pr}_1 :\Sigma_{\bf g}\times [-1,1]\to \Sigma_{\bf g}$ is now
 holomorphic,
with rational curves as fibers. To get a minimal model for $M$, we can
proceed by observing that, for any ``puncture point'' $q_j$,
the inverse image $\pi^{-1}(\{ \mbox{pr}_1(q_j)\}\times
[-1, t(q_j)])$
 of the  vertical line segment  joining the lower boundary of
$\Sigma_{\bf g}\times [-1,1]$
to  $q_j$ is a rational curve  in $M$,
and, provided the segment
does not pass through any other puncture point, this rational curve
is smooth, with self-intersection $-1$. (In the non-generic situation in which
 several
of the puncture points project to the same point of $\Sigma_{\bf g}$,
the line segment between any
two such consecutive points similarly
corresponds to a smooth rational curve in $M$ of
self-intersection $-2$.) By blowing down all such $(-1)$-curves
(and  iteratively blowing down the  $(-1)$-curves that  then arise from
$(-2)$-curves in the non-generic case) we eventually
arrive at a minimal model ${\Bbb P}({\cal L}\oplus {\cal O})\to \Sigma_{\bf g}$
with holomorphic vector field, where all the blow-ups occur at the
the zero section of ${\cal L}$. The proper transform $C_0$
of the zero section now  corresponds to $t=-1$, whereas the infinity section
$C_{\infty }$ corresponds to $t=+1$. Meanwhile,
the line bundle ${\cal L}^{\ast}$
is  exactly the holomorphic line bundle  associated to the
$U(1)$-connection obtained by restricting  $(M_0, \theta)$ to $\Sigma_{\bf g}
\times \{ 1-\epsilon\}$, so that equation (\ref{int}) yields
\be \mbox{deg} ({\cal L})=\sum_{j=1}^mw_j~ .\label{deg}\ee
But a different choice of $M_0$ would change this
$U(1)$-connection by twisting it by an arbitrary flat $U(1)$-connection;
since $\mbox{Pic}_0(\Sigma_{\bf g})=H^1(\Sigma_{\bf g},
{\cal O})/H^1(\Sigma_{\bf g}, {\Bbb Z})$ is canonically
 identified with $H^1(\Sigma_{\bf g}, S^1)=
H^1(\Sigma_{\bf g}, {\Bbb R})/H^1(\Sigma_{\bf g}, {\Bbb Z})$
by the Hodge decomposition,
this means that the holomorphic line-bundles which arise
 for a fixed configuration $q_1, \ldots , q_m$ fill out the
 entire connected component of $\mbox{Pic}(\Sigma_{\bf g})$
specified by the degree formula (\ref{deg}).
And since the fiber-wise K\"ahler form is just\footnote{This is again
a manifestation of  the ``Archimedes Principle'' \cite{archie}\cite{atiyah}
for symplectic torus actions.}
$$\left.\omega \right|_{ fiber}= dt\wedge\theta ~,$$
the area of the holomorphic curve $E_j=\pi^{-1}(\{ \mbox{pr}_1 (q_j)\}\times
[-1, t(q_j)])$ is just
$$\int_{E_j}\omega = 2\pi (t(q_j)-(-1))= 4\pi w_j~ .$$
Since, by the same reasoning,
 the typical fiber $F$ of $M\to \Sigma_{\bf g}$ has area
$4\pi$, the numbers $w_j$ are precisely the weights we associated with the
exceptional divisors in \S \ref{foo}, and equation (\ref{deg}) is therefore,
by Proposition \ref{foot}, just
the assertion that the Futaki invariant vanishes--- as of course it must,
since our K\"ahler manifold has constant scalar curvature zero! Since
we are free to choose the numbers $w_j$,  subject only to the constraint
(\ref{foutre}), and since, by multiplying $g$ by an arbitrary constant,
we can make the typical fiber $F$ have any area we choose,
the above explicit construction produces  a scalar-flat
K\"ahler metric in any K\"ahler class on $M$ such that both
$c_1\cdot [\omega ]$ and  the Futaki invariant are zero:

\begin{thm} Let $M$ be a compact complex surface with $a(M)\neq 0$.
Then a K\"ahler class $[\omega ]\in
H^{1,1}(M)$ contains a scalar-flat K\"ahler metric iff  the
total scalar curvature, the Matsushima-Lichnerowicz obstruction,
and the Futaki invariant  all vanish.
When such a metric exists, it is unique modulo
biholomorphisms of $M$. \label{class}
\end{thm}
\begin{proof} For the existence part, it remains only to observe that
the only  non-ruled cases are tori and   hyperelliptic  surfaces,
and these admit flat
metrics in every K\"ahler class.
For the uniqueness result, which we shall never use in this article,
we refer the reader to  \cite{L5}, Theorem 3.
\end{proof}

\noindent
The following simple application will prove to be particularly useful:
\begin{cor}
Let the product surface $\Sigma_{\bf g}\times {\Bbb CP}_1$ be
blown up at any $k>0$ points along $\Sigma_{\bf g}\times\{ [1:0]\}$
and any
$\ell > 0$ points along $\Sigma_{\bf g}\times\{ [0:1]\}$. (Some or all of the
given  points are allowed to coincide, but in this case the iterated blow-ups
are required to occur along proper transforms of the
$\Sigma_{\bf g}$ or ${\Bbb CP}_1$ factors.)
The resulting surface then admits scalar-flat K\"ahler metrics.
\label{triv}
\end{cor}
\begin{proof} If the set of blown-up points  is $$\{ (p_1, [1:0]),
\ldots , (p_k, [1:0]), (q_1, [0:1]), \ldots , (q_{\ell}, [0:1])\}~,$$
then, letting ${\cal L}\to \Sigma_{\bf g}$ denote the divisor
line bundle of $\{ q_1, \ldots , q_{\ell}\}$,
the surface in question can also be described as the blow-up
of ${\Bbb P}(\O \oplus {\cal L})$ at the points $\{ p_1,
\ldots , p_k,q_1, \ldots , q_{\ell}\}$ on the zero section.
Since $0<\mbox{deg}({\cal L})=\ell<m=\ell +k$, the result follows from
Corollary \ref{tri} and
Theorem \ref{class}.
\end{proof}

\noindent
In light of \cite{burnsbart}, the following  restatement of the
above theorem seems particularly tantalizing:
\begin{cor} An admissible  K\"ahler class on a (blown-up) ruled surface
with periodic
holomorphic vector field contains a scalar-flat K\"ahler metric iff
the corresponding parabolic bundle is quasi-stable.
\label{corn}
\end{cor}

\subsection{Scalar-Flat Metrics on Generic Ruled Surfaces }\label{tag}

The key technical result of this paper is as follows:
\begin{thm} Let $\varpi : {\cal M}\to {\cal U}$ be a family
of non-minimal ruled surfaces of genus $\geq 2$. Suppose that,
for some $o\in {\cal U}$, the corresponding fiber $M=M_o:= \varpi^{-1}(o)$
admits a  scalar-flat K\"ahler metric. Then there is a neighborhood
$\tilde{\cal U}\subset {\cal U}$ of  $o$ such that $M_t:= \varpi^{-1}(t)$
admits a scalar-flat K\"ahler metric for all $t\in \tilde{\cal U}$. Moreover,
relative to  local trivializations of the real-analytic fiber-bundle
underlying $\varpi$,
these metrics can be chosen so as to depend real-analytically on $t$.
\label{tech}
\end{thm}
\begin{proof} Let $Z$ be the twistor space of $M$, $D=M\coprod \bar{M}$
its standard divisor,
and $\sigma$ its real structure.
Using Theorem \ref{dt5} and Corollary \ref{van}, we have
$H^{2}(Z,{\Theta}_{Z}\otimes {\cal I}_D)=H^{2}(Z, {\Theta}_{Z,D})=0$.
The first statement can be reinterpreted on $Z/\sigma$ as saying that
$H^{2}(Z/\sigma,\Re{\Theta}_{Z}\otimes {\cal I}_D)=0$, so that
the long exact sequence induced by the short exact sequence
$$0\to \Re [{\Theta}_{P}\otimes {\cal I}_D ]
 \to \Re{\Theta}_{P}\to \Theta_M\to 0$$
on $P=Z/\sigma$ predicts that the natural restriction map
$$H^1(Z/\sigma , \Re{\Theta})\to H^1(M, \Theta_M)$$
is surjective. Applying this morphism to restrict the
Kodaira-Spencer map to $M$, we see that
 the versal family for $(Z/\sigma , M)$
given  by Theorem \ref{ks3}  induces a complete deformation of $M$--- i.e.
a deformation of $M$ which
contains a versal deformation as a subspace.
Thus any small deformation of $M$ can be extended as a deformation
of $(Z, D, \sigma)$. Applying  Theorem \ref{Pen} then finishes the proof.
\end{proof}

\begin{remark} The analogous statement fails \cite{burnsbart}
for {\em minimal}
ruled surfaces.
\end{remark}

\noindent We now prove our main result:
\begin{thm} Let $M$ be any ruled surface of genus $\geq 2$.
If $M$ is blown up at sufficiently
many points, the resulting complex surface $\tilde{M}$
admits scalar-flat K\"ahler metrics.
\label{mainline}
\end{thm}
\begin{proof} Since any ruled surface \cite{bpv} is bimeromorphic to a
product surface, some blow-up $\tilde{M}$ of the given surface
 $M$ is biregularly
equivalent to an iterated blow-up of $\Sigma_{\bf g}\times {\Bbb CP}_1$,
${\bf g}\geq 2$,
at  $r_1\ldots r_m\in \Sigma_{\bf g}\times {\Bbb CP}_1$, where
repetition of a point  indicates that we are  also
 given certain  directional information at the multiple point.
By blowing up extra points if necessary, we may assume that
the  projection of   $\{r_1\ldots r_m \}$ to
${\Bbb CP}_1$ consists of more than one point. By changing to another
homogeneous coordinate system $[\zeta_1 : \zeta_2]$
on  ${\Bbb CP}_1$, we may also assume that $[1:0]$ is in the image
of  $\{r_1\ldots r_m \}$. For $t\in {\Bbb C}$, and if $r_j$
 projects to $[1:0]$, set $r_j(t):= r_j$; otherwise
let $r_j(t):= \mu_t (r_j)$, where
\bea \mu_t : \Sigma_{\bf g}\times ({\Bbb CP}_1- \{ [1:0]\}  )&\to&
 \Sigma_{\bf g}\times {\Bbb CP}_1
\\(p , [\zeta_1 : \zeta_2])&\mapsto&
(p , [t\zeta_1 : \zeta_2])~ .
\eea
Define a family $\varpi :{\cal M}\to {\Bbb C}$ of complex surfaces
by blowing up $\Sigma_{\bf g}\times {\Bbb CP}_1
\times {\Bbb C}$ along the graphs of $t\mapsto r_j(t) $.
For $t\neq 0$, the manifold $M_t$ is then biholomorphic to $\tilde{M}$.
However, for $t=0$, the fiber is a ruled surface with a holomorphic
vector field, namely the blow-up of $\Sigma_{\bf g}\times {\Bbb CP}_1$
at  non-empty collections of points along $\Sigma_{\bf g}\times \{ [0:1]\}$
and $\Sigma_{\bf g}\times \{ [1:0]\}$. (When several of the blown-up points
project to the same point of $\Sigma_{\bf g}$, one should think
of the blow-ups as happening successively rather than simultaneously;
for $t=0$, we are, at
each stage,  blowing up the previous surface
at a zero of the vector field $\Xi = \zeta_1\partial /
\partial \zeta_1$, and $\Xi$ therefore lifts to the blow-up.)
Corollary \ref{triv} then asserts that $M_0$ admits scalar-flat K\"ahler
metrics. The result therefore follows immediately from
Theorem \ref{tech}.

\end{proof}
This can be repackaged as follows:
\begin{cor} {\rm \bf (Main Theorem)}
Let $M$ be a compact complex surface which admits
a K\"ahler metric whose scalar curvature has integral zero. Suppose
$\pi_1 (M)$ does not contain an Abelian subgroup of finite index. Then
if $M$ is blown up at sufficiently many points, the
resulting surface $\tilde{M}$ admits scalar-flat K\"ahler metrics.
\end{cor}
\begin{proof} By Theorem \ref{sd3}, a complex surface with admissible
K\"ahler class is either ruled or covered by a torus or K3. Thus the
fundamental group hypothesis forces $M$ to be a ruled surface of
genus $\geq 2$. The statement therefore follows from Theorem
\ref{mainline}.
\end{proof}

We now conclude this article
with a pair of conjectures intended to remind the reader that
the study of scalar-flat K\"ahler surfaces
 is still in its infancy.  First, one would like to understand what happens
for ruled surfaces of  genus 0 and 1. It is our hope and expectation
that the genus hypothesis in Theorem \ref{mainline} is actually
superfluous:

\begin{conjecture} The blow-up of {\em any} ruled surface at sufficiently
many points admits scalar-flat K\"ahler metrics. \label{better}
\end{conjecture}

On the other hand, one would really like to understand {\em precisely} when an
admissible class contains a scalar-flat K\"ahler metric.
In light of Corollary \ref{corn} and
known results on relatively  minimal  ruled surfaces \cite{burnsbart},
 the following would seem very natural:

\begin{conjecture} An admissible  K\"ahler class on a (blown-up) ruled surface
of genus $\geq 2$ contains a scalar-flat K\"ahler metric iff
the corresponding parabolic bundle is quasi-stable.
\end{conjecture}

\vfill

\noindent {\bf Acknowledgements.}  The authors would like to thank the
Australian Research Council for funding their visits to the
University of Adelaide, where much of the actual writing was done.
They would also like to thank
N P.\ Buchdahl, M.G.\ Eastwood, M.S.\ Narasimhan,
Y.-T.\ Siu, and G.\ Tian for their helpful comments and
suggestions.

\pagebreak

\end{document}